\documentclass[a4paper,11pt]{amsart}

\usepackage{graphicx}
\usepackage{mathptmx}
\usepackage{amsmath}
\usepackage{amssymb}
\usepackage{enumitem}
\usepackage{xcolor}

\newmuskip\pFqmuskip

\newcommand*\pFq[6][8]{%
  \begingroup % only local assignments
  \pFqmuskip=#1mu\relax
  \mathcode`=\string"8000
  \begingroup\lccode`\~=`\,
  \lowercase{\endgroup\let~}\pFqcomma
  F^{#2}_{#3}{\left(\genfrac..{0pt}{}{#4}{#5}\bigg|#6\right)}%
  \endgroup
}
\newcommand{\pFqcomma}{\mskip\pFqmuskip}

\numberwithin{equation}{section}
\newtheorem{theorem}{Theorem}[section]
\newtheorem{proposition}[theorem]{Proposition}

\newtheorem{remark}[theorem]{Remark}
\newtheorem{lemma}[theorem]{Lemma}

\begin{document}

\pagenumbering{arabic}
\pagestyle{headings}

\newcommand{\DPB}[4]{P\beta_{#1}^{(#2)}(#3,#4)}

\title{Representations by probabilistic Bernoulli and degenerate Bernoulli polynomials}
\author{Dae San Kim}
\address{Department of Mathematics, Sogang University, Seoul 121-742, Republic of Korea}
\email{dskim@sogang.ac.kr}

\author{Taekyun Kim}
\address{Department of Mathematics, Kwangwoon University, Seoul 139-701, Republic of Korea}
\email{tkkim@kw.ac.kr}


\begin{abstract}
We investigate the representation of arbitrary polynomials using probabilistic Bernoulli and degenerate Bernoulli polynomials associated with a random variable $Y$, whose moment generating function exists in a neighborhood of the origin. In addition, this paper explores the problem of representing arbitrary polynomials in terms of their higher-order counterparts. We develop explicit formulas for those representations with the help of umbral calculus and illustrate our results for several discrete and continuous random variables $Y$.
\end{abstract}

\maketitle

\markboth{\centerline{\scriptsize Representations by probabilistic Bernoulli and degenerate Bernoulli polynomials}}
{\centerline{\scriptsize Dae San Kim and Taekyun Kim}}

\section{Introduction and preliminaries}
Let $Y$ be a random variable with a moment generating function existing in a neighborhood of the origin (see \eqref{1a}). The study of degenerate versions of special polynomials and numbers, originating with Carlitz's work on degenerate Bernoulli and Euler polynomials [4], has seen a resurgence of interest [9,13,19,20,23, and references therein]. Similarly, probabilistic extensions of special polynomials and numbers have been extensively researched [1-3,6,18,21,22,32]. \par
This paper investigates expressing arbitrary polynomials as linear combinations of probabilistic Bernoulli polynomials associated with $Y$, $B_{n}^{Y}(x)$, and probabilistic degenerate Bernoulli polynomials associated with $Y$, $\beta_{n}^{Y}(x)$. We also address representing polynomials in terms of their higher-order counterparts. Utilizing umbral calculus [8,27,28,31], we derive formulas for the coefficients in these linear combinations. As examples, we express $x^{n}$ as linear combinations of $B_{k}^{Y}(x)$ and $\beta_{k}^{Y}(x)$, for various discrete and continuous random variables $Y$. This requires explicit expressions for the probabilistic Stirling numbers of the first kind associated with $Y$, $S_{1}^{Y}(n,k)$, and the probabilistic degenerate Stirling numbers of the first kind associated with $Y$, $S_{1,\lambda}^{Y}(n,k)$ [14]. Then, as illustrations for our results, we express $x^{n}$ as linear combinations of $B_{k}^{Y}(x)$ and $\beta_{k}^{Y}(x)$, for several discrete and continuous random variables $Y$. Crucially, $S_{1}^{Y}(n,k)$ and $S_{2}^{Y}(n,k)$, and $S_{1,\lambda}^{Y}(n,k)$ and $S_{2,\lambda}^{Y}(n,k)$, satisfy orthogonality and inverse relations (Propositions 1.1 and 1.2), which are, as inversions are needed, essential for our problems. In contrast, the definitions of $S_{1}^{Y}(n,k)$ in [2] and $S_{1,\lambda}^{Y}(n,k)$ in [18], based on cumulant generating functions, respectively together with $S_{2}^{Y}(n,k)$ and $S_{2,\lambda}^{Y}(n,k)$, do not possess these properties. To illustrate our results, we express $x^{n}$ as linear combinations of $B_{k}^{Y}(x)$ and of $\beta_{k,\lambda}^{Y}(x)$ for several random variables $Y$, relying on explicit computations of $S_{1}^{Y}(n,k)$ and $S_{1,\lambda}^{Y}(n,k)$ [14]. Consider, for example, the exponential random variable $Y$ with parameter $\alpha >0$, whose probability density function is given by [29]
\begin{equation*}
f(y)=\left\{\begin{array}{ccc}
\alpha e^{-\alpha y}, & \textrm{if \,\,$y \ge 0$,} \\
0, & \textrm{if\,\, $y<0$}.
\end{array}\right.
\end{equation*}
Then our result in this case is as follows:
\begin{align*}
x^{n}&= \frac{1}{\alpha}B_{0}^{Y}(x)+\sum_{k=1}^{n}\Big\{\frac{1}{k}\sum_{j=k-1}^{n-1}(-1)^{j-k+1}\binom{j}{k-1} \\
&\quad\quad \times (j-1)_{j-k+1}\alpha^{k-1}\frac{1}{j!}\Delta^{j+1}0^{n}\Big\}B_{k}^{Y}(x), \\
x^{n}&=\frac{1}{\alpha}\sum_{r=0}^{n}\sum_{l=0}^{r}\binom{n}{r}S_{2}(r,l)(-1)^{l-r}(\alpha \lambda)^{l} B_{l}\beta_{0,\lambda}^{Y}(x) \\
&\quad +\sum_{k=1}^{n}\Big\{\frac{1}{k}\sum_{j=k-1}^{n-1}\sum_{m=k-1}^{j}\binom{j}{m}(-1)^{j-m}(j-1)_{j-m}  \\
& \quad \quad \quad \times \alpha^{m}\lambda^{m-k+1}S_{2}(m,k-1)\frac{1}{j!}\Delta^{j+1}0^{n} \Big\}\beta_{k,\lambda}^{Y}(x).
\end{align*} \par
Here we mention some of the previous results. Applying the formula in \eqref{21-1c} to $p(x)=\sum_{k=1}^{n-1}\frac{1}{k(n-k)}B_{k}(x)B_{n-k}(x)$ yields the following identity [17, Theorem 4]:
\begin{equation}
\sum_{k=1}^{n-1}\frac{B_{k}(x)B_{n-k}(x)}{k(n-k)}=\frac{2}{n}\sum_{k=0}^{n-2}\frac{1}{n-k}\binom{n}{k}B_{n-k}B_{k}(x)+\frac{2}{n}H_{n-1}B_{n}(x), \label{-1a}
\end{equation}
where $n \ge 2$, and $H_{n}=1+\frac{1}{2}+ \cdots +\frac{1}{n}$. Substituting
$x=0$ and $x=\frac{1}{2}$ into equation \eqref{-1a} yields, respectively, a variant of Miki's identity and the Faber-Pandharipande-Zagier (FPZ) identity (see [16]). Notably, our derivation of these identities relies on the straightforward formula \eqref{21-1c}, which involves only derivatives and integrals of the given polynomials. This contrasts sharply with existing proofs, which are considerably more complex. For Miki's identity, Miki [24] employed a formula for the Fermat quotient $\frac{a^{p}-a}{p}$ modulo $p^{2}$, Shiratani and Yokoyama [30] utilized $p$-adic analysis, and Gessel [12] leveraged two distinct expressions for Stirling numbers of the second kind, $S_{2}\left(n,k\right)$. Similarly, Dunne and Schubert [10] derived the FPZ identity using asymptotic expansions of special polynomials arising from quantum field theory. Zagier also provided a proof in the appendix of [11]. It's worth noting that Faber and Pandharipande initially conjectured relations between Hodge integrals in Gromov-Witten theory, which necessitated the FPZ identity, in 1998. Further related work can be found in [5,15,17]. The representation of arbitrary polynomials in terms of $\beta_{k,\lambda}(x)$ using formulas \eqref{34c} and \eqref{35c} is discussed in [16]. Specifically, applying these formulas to $p(x)=\sum_{k=1}^{n-1}\frac{1}{k(n-k)}B_{k}(x)B_{n-k}(x)$ leads to the following identity (see [16, Examples (b)]):
\begin{align*}
&\sum_{k=1}^{n-1}\frac{1}{k(n-k)}B_{k}(x)B_{n-k}(x)\\
&=\frac{2}{n}\bigg\{\sum_{l=0}^{n-2}\frac{1}{n-l}\binom{n}{l}B_{n-l}\lambda^{l}B_{l}+H_{n-1}\lambda^{n}B_{n}\bigg\}\beta_{0,\lambda}(x)\\
&\quad\quad+\frac{2}{n}\sum_{k=1}^{n}\frac{1}{k!\lambda^{k-1}}\bigg\{\sum_{l=1}^{n-2}\frac{l}{n-l}\binom{n}{l}B_{n-l}\Delta_{\lambda}^{k-1}0^{l-1}+n H_{n-1}\Delta_{\lambda}^{k-1}0^{n-1}\bigg\}\beta_{k,\lambda}(x),
\end{align*}
where  $n \ge 2$. \par
This paper is organized as follows. Section 1 reviews the necessary preliminaries. Section 2 provides a succinct overview of umbral calculus. Our main contributions are presented in Sections 3 and 4, where we derive novel formulas expressing arbitrary polynomials in terms of probabilistic and higher-order probabilistic Bernoulli polynomials, and their degenerate counterparts, associated with a random variable $Y$. Section 5 illustrates these results with concrete examples, specifically expressing $x^{n}$ in terms of $B_{n}^{Y}(x)$ and $\beta_{n,\lambda}^{Y}(x)$,
for various choices of $Y$. Finally, Section 6 summarizes our findings. General background references include [7,25,26,29]. The rest of Section 1 details the required background information. \par
\vspace{0.1in}
Let $Y$ be a random variable whose moment generating function exists in a neighborhood of the origin:
\begin{equation}
E[e^{Yt}]=\sum_{n=0}^{\infty}E[Y^{n}]\frac{t^{n}}{n!}\,\,\, \mathrm{exists,\,\, for}\,\,|x|<r, \label{1a}
\end{equation}
for some positive real number $r$. \\
Let $(Y_{j})_{j \ge 1}$ be a sequence of mutually independent copies of the random variable $Y$, and let
\begin{equation}
S_{k}=Y_{1}+Y_{2}+\cdots +Y_{k}, \quad (k \ge 1), \quad  S_{0}=0. \label{2a}
\end{equation}
The probabilistic Stirling numbers of the second kind associated with $Y$, $S_{2}^{Y}(n,k)$, are given by (see \eqref{2a})
\begin{align}
&\frac{1}{k!}(E[e^{Yt}]-1)^{k}=\sum_{n=k}^{\infty}S_{2}^{Y}(n,k)\frac{t^{n}}{n!}, \label{3a}\\
&S_{2}^{Y}(n,k)=\frac{1}{k!}\sum_{j=0}^{k}\binom{k}{j}(-1)^{k-j}E[S_{j}^{n}]. \nonumber
\end{align}
From the definition in \eqref{3a}, it is immediate to see that
\begin{equation}
S_{2}^{Y}(k,k)=E[Y]^{k}. \label{4a}
\end{equation} \par
Assume from now on that
\begin{equation}
E[Y] \ne 0. \label{5a}
\end{equation}
We introduce the notation:
\begin{equation}
e_{Y}(t)=E[e^{Yt}]-1. \label{6a}
\end{equation}
Then we have
\begin{equation}
\frac{1}{k!}\big(e_{Y}(t)\big)^{k}=\sum_{n=k}^{\infty}S_{2}^{Y}(n,k)\frac{t^{n}}{n!}. \label{7a}
\end{equation}
If $f(t)=\sum_{n=0}^{\infty}a_{n}\frac{t^{n}}{n!}$ is a delta series, namely $a_{0}=0$ and $a_{1} \ne 0$, then the compositional inverse $\bar{f}(t)$ of $f(t)$ satisfying $f(\bar{f}(t))=\bar{f}(f(t))=t$ exists. Note that, as $e_{Y}(t)=E[Y]t+\sum_{m=2}^{\infty}E[Y^{m}]\frac{t^{m}}{m!}$ and $E[Y] \ne 0$ (see \eqref{5a}, \eqref{6a}), $e_{Y}(t)$ is a delta series. \par
Now, we define the  probabilistic Stirling numbers of the first kind associated with $Y$ by: for $k \ge 0$,
\begin{equation}
\frac{1}{k!}\big(\bar{e}_{Y}(t)\big)^{k}=\sum_{n=k}^{\infty}S_{1}^{Y}(n,k)\frac{t^{n}}{n!}, \label{8a}
\end{equation}
where $\bar{e}_{Y}(t)$ is the compositional inverse of $e_{Y}(t)$. \\
In addition, as usual, we agree that
\begin{equation}
S_{2}^{Y}(n,k)=S_{1}^{Y}(n,k)=0,\,\, \mathrm{if}\,\, k>n\,\,\mathrm{or}\,\, k <0. \label{9a}
\end{equation}
Note that $S_{2}^{Y}(n,k)=S_{2}(n,k),\,\, S_{1}^{Y}(n,k)=S_{1}(n,k)$, for $Y=1$. \\
Here $S_{2}(n,k)$ are the Stirling numbers of the second kind defined by
\begin{align}
&x^{n}=\sum_{k=0}^{n}S_{2}(n,k)(x)_{k},\label{10a} \\
&\frac{1}{k!}(e^{t}-1)^{k}=\sum_{n=k}^{\infty}S_{2}(n,k)\frac{t^{n}}{n!}, \nonumber
\end{align}
and $S_{1}(n,k)$ are the Stirling numbers of the first kind defined as
\begin{align}
&(x)_{n}=\sum_{k=0}^{n}S_{1}(n,k)x^{k}, \label{11a} \\
&\frac{1}{k!}\big(\log(1+t)\big)^{k}=\sum_{n=k}^{\infty}S_{1}(n,k)\frac{t^{n}}{n!}, \nonumber
\end{align}
where $(x)_{n}$ are the falling factorials given by
\begin{equation}
(x)_{0}=1, \quad (x)_{n}=x(x-1)\cdots(x-n+1),\quad (n \ge 1). \label{12a}
\end{equation} \par
Using the definitions in \eqref{7a} and \eqref{8a}, one shows that $S_{2}^{Y}(n,k)$ and $S_{1}^{Y}(n,k)$ satisfy the orthogonality relations in (a) of Proposition 1.1, from which the inverse relations in (b) and (c) follow.
\begin{proposition}
The following orthogonality and inverse relations are valid for $S_{1}^{Y}(n,k)$ and $S_{2}^{Y}(n,k)$.
\begin{flalign*}
&(a)\,\, \,\sum_{k=l}^{n} S_{2}^{Y}(n,k)S_{1}^{Y}(k,l)=\delta_{n,l}, \quad \sum_{k=l}^{n} S_{1}^{Y}(n,k)S_{2}^{Y}(k,l)=\delta_{n,l}, \\
&(b)\,\, a_{n}=\sum_{k=0}^{n}S_{2}^{Y}(n,k) b_{k}\,\, \iff \,\, b_{n}=\sum_{k=0}^{n}S_{1}^{Y}(n,k)a_{k}, \\
&(c)\,\, a_{n}=\sum_{k=n}^{m}S_{2}^{Y}(k,n)b_{k} \,\, \iff \,\, b_{n}=\sum_{k=n}^{m}S_{1}^{Y}(k,n)a_{k}. &&
\end{flalign*}
\end{proposition}
Let $\lambda$ be any nonzero real number. Then $e_{\lambda}^{x}(t)$ are the degenerate exponentials defined by
\begin{equation}
e_{\lambda}^{x}(t)=(1+\lambda t)^{\frac{x}{\lambda}}=\sum_{n=0}^{\infty}(x)_{n,\lambda}\frac{t^{n}}{n!}, \quad e_{\lambda}(t)=e_{\lambda}^{1}(t),\quad (\mathrm{see}\ [18,21]), \label{13a}
\end{equation}
where $(x)_{n,\lambda}$ are the degenerate falling factorials given by
\begin{equation}
(x)_{0,\lambda}=1,\ (x)_{n,\lambda}=x(x-\lambda)\cdots(x-(n-1)\lambda),\ (n\ge 1).\label{14a}
\end{equation}
Note here that $\lim_{\lambda\rightarrow 0}e_{\lambda}^{x}(t)=e^{xt}$. \\
The probabilistic degenerate Stirling numbers of the second kind associated $Y$, $S_{2, \lambda}^{Y}(n,k)$, are defined by
\begin{align}
&\frac{1}{k!}\big(E[e_{\lambda}^{Y}(t)]-1\big)^{k}=\sum_{n=k}^{\infty}S_{2,\lambda}^{Y}(n,k)\frac{t^{n}}{n!}, \label{15a} \\
&S_{2,\lambda}^{Y}(n,k)=\frac{1}{k!}\sum_{j=0}^{k}\binom{k}{j}(-1)^{k-j}E[(S_{j})_{n,\lambda}]. \nonumber
\end{align} \par
To define the probabilistic degenerate Stirling numbers of the first kind associated with $Y$, we let
\begin{equation}
e_{Y,\lambda}(t)=E[e_{\lambda}^{Y}(t)]-1. \label{16a}
\end{equation}
Then we have
\begin{equation}
\frac{1}{k!}(e_{Y,\lambda}(t))^{k}=\sum_{n=k}^{\infty}S_{2,\lambda}^{Y}(n,k)\frac{t^{n}}{n!}. \label{17a}
\end{equation}
Noting that $e_{Y,\lambda}(t)=E[Y]t+\sum_{m=2}^{\infty}E[(Y)_{m,\lambda}]\frac{t^{m}}{m!}$ is a delta series (see \eqref{14a}), we define the probabilistic degenerate Stirling numbers of the first kind associated with $Y$, $S_{1,\lambda}^{Y}(n,k)$, by
\begin{equation}
\frac{1}{k!}\big(\bar{e}_{Y,\lambda}(t)\big)^{k}=\sum_{n=k}^{\infty}S_{1,\lambda}^{Y}(n,k)\frac{t^{n}}{n!}, \label{18a}
\end{equation}
where $\bar{e}_{Y,\lambda}(t)$ is the compositional inverse of $e_{Y,\lambda}(t)$. \\
Note that $S_{2,\lambda}^{Y}(n,k)=S_{2,\lambda}(n,k)$ and $S_{1,\lambda}^{Y}(n,k)=S_{1,\lambda}(n,k)$, for $Y=1$.
Here $S_{2,\lambda}(n,k)$ are the degenerate Stirling numbers of the second kind defined by
\begin{align}
&(x)_{n,\lambda}=\sum_{k=0}^{n}S_{2,\lambda}(n,k)(x)_{k}, \label{19a} \\
&\frac{1}{k!}(e_{\lambda}(t)-1)^{k}=\sum_{n=k}^{\infty}S_{2,\lambda}(n,k)\frac{t^{n}}{n!}, \nonumber
\end{align}
and $S_{1,\lambda}(n,k)$ are the degenerate Stirling numbers of the first kind given by
\begin{align}
&(x)_{n}=\sum_{k=0}^{n}S_{1,\lambda}(n,k)(x)_{k,\lambda}, \label{20a} \\
&\frac{1}{k!}\big(\log_{\lambda}(1+t)\big)^{k}=\sum_{n=k}^{\infty}S_{1,\lambda}(n,k)\frac{t^{n}}{n!}, \nonumber
\end{align}
where $\log_{\lambda}(t)$ are the degenerate logarithm defined by
\begin{equation}
\log_{\lambda}(t)=\frac{1}{\lambda}\big(t^{\lambda}-1\big). \label{21a}
\end{equation}
Note here that the degenerate exponential $e_{\lambda}(t)$ in \eqref{13a} and the degenerate logarithm $\log_{\lambda}(t)$ in \eqref{21a} are compositional inverses to each other so that
\begin{equation}
e_{\lambda}\big(\log_{\lambda}(t) \big)=\log_{\lambda}\big(e_{\lambda}(t) \big)=t. \label{22a}
\end{equation} \par
From \eqref{17a} and \eqref{18a}, one shows that $S_{2,\lambda}^{Y}(n,k)$ and $S_{1,\lambda}^{Y}(n,k)$ satisfy the orthogonality relations in (a) of Proposition 1.2, from which the inverse relations in (b) and (c) follow.
\begin{proposition}
The following orthogonality and inverse relations are valid for $S_{1,\lambda}^{Y}(n,k)$ and $S_{2,\lambda}^{Y}(n,k)$.
\begin{flalign*}
&(a)\,\, \sum_{k=l}^{n} S_{2,\lambda}^{Y}(n,k)S_{1,\lambda}^{Y}(k,l)=\delta_{n,l}, \quad \sum_{k=l}^{n} S_{1,\lambda}^{Y}(n,k)S_{2,\lambda}^{Y}(k,l)=\delta_{n,l},\\
&(b)\,\, a_{n}=\sum_{k=0}^{n}S_{2,\lambda}^{Y}(n,k) b_{k}\,\, \iff \,\, b_{n}=\sum_{k=0}^{n}S_{1,\lambda}^{Y}(n,k)a_{k}, \\
&(c)\,\, a_{n}=\sum_{k=n}^{m}S_{2,\lambda}^{Y}(k,n)b_{k} \,\, \iff \,\, b_{n}=\sum_{k=n}^{m}S_{1,\lambda}^{Y}(k,n)a_{k}. &&
\end{flalign*}
\end{proposition}
\noindent Notice that the Stirling numbers of both kinds and the degenerate Stirling numbers of both kinds satisfy orthogonality and inversion relations (see Propositions 1.1 and 1.2 with $Y=1$). \par
The probabilistic Bernoulli polynomials associated with $Y$, $B_{n}^{Y}(x)$, are defined by
\begin{equation}
\frac{t}{E[e^{Yt}]-1}\big(E[e^{Yt}] \big)^{x}=\sum_{n=0}^{\infty}B_{n}^{Y}(x)\frac{t^{n}}{n!}. \label{23a}
\end{equation}
More generally, for any nonnegative integer $r$, the probabilistic Bernoulli polynomials of order $r$ associated with $Y$, $B_{n}^{Y,(r)}(x)$, are given by
\begin{equation}
\bigg(\frac{t}{E[e^{Yt}]-1}\bigg)^{r}\big(E[e^{Yt}] \big)^{x}=\sum_{n=0}^{\infty}B_{n}^{Y,(r)}(x)\frac{t^{n}}{n!}. \label{24a}
\end{equation}
The probabilistic degenerate Bernoulli polynomials associated with $Y$, $\beta_{n,\lambda}^{Y}(x)$, are defined by
\begin{equation}
\frac{t}{E[e_{\lambda}^{Y}(t)]-1}\big(E[e_{\lambda}^{Y}(t)] \big)^{x}=\sum_{n=0}^{\infty}\beta_{n,\lambda}^{Y}(x)\frac{t^{n}}{n!}. \label{25a}
\end{equation}
More generally, for any nonnegative integer $r$, the probabilistic degenerate Bernoulli polynomials of order $r$ associated with $Y$, $\beta_{n,\lambda}^{Y,(r)}(x)$, are given by
\begin{equation}
\bigg(\frac{t}{E[e_{\lambda}^{Y}(t)]-1}\bigg)^{r}\big(E[e_{\lambda}^{Y}(t)] \big)^{x}=\sum_{n=0}^{\infty}\beta_{n,\lambda}^{Y,(r)}(x)\frac{t^{n}}{n!}. \label{26a}
\end{equation} \par
We recall some notations and facts about forward differences. Let $f$ be any complex-valued function of the real variable $x$. Then, for any real number $a$, the forward difference $\Delta_{a}$ is given by
\begin{equation}
\Delta_{a}f(x)=f(x+a)-f(x).\label{27a}
\end{equation}
If $a=1$, then we let
\begin{equation}
\Delta f(x)=\Delta_{1}f(x)=f(x+1)-f(x).\label{28a}
\end{equation}
In general, the $n$th order forward differences are given by
\begin{equation}
\Delta_{a}^{n}f(x)=\sum_{i=0}^{n}\binom{n}{i} (-1)^{n-i}f(x+ia).\label{29a}
\end{equation}
For $a=1$, we have
\begin{equation}
\Delta^{n}f(x)=\sum_{i=0}^{n}\binom{n}{i} (-1)^{n-i}f(x+i).\label{30a}
\end{equation} \par
For $Y=1$, $B_{n}(x)=B_{n}^{1}(x)$ are the Bernoulli polynomials defined by
\begin{equation}
\frac{t}{e^t-1}e^{xt}=\sum_{n=0}^{\infty}B_{n}(x)\frac{t^n}{n!}.\label{31a}
\end{equation}
When $x=0$, $B_n=B_n(0)$ are called the Bernoulli numbers. \\
For $Y=1$, $B_{n}^{(r)}(x)=B_{n}^{1,(r)}(x)$ are the Bernoulli polynomials of order $r$,
given by
\begin{equation}
\bigg(\frac{t}{e^t-1}\bigg)^{r}e^{xt}=\sum_{n=0}^{\infty}B_{n}^{(r)}(x)\frac{t^n}{n!}.\label{32a}
\end{equation}
For $Y=1$, $\beta_{n,\lambda}(x)=\beta_{n,\lambda}^{1}(x)$ are the degenerate Bernoulli polynomials given by
\begin{equation}
\frac{t}{e_{\lambda}(t)-1}e_{\lambda}^{x}(t)=\sum_{n=0}^{\infty}\beta_{n,\lambda}(x)\frac{t^n}{n!}.\label{33a}
\end{equation}
When $x=0$, $\beta_{n,\lambda}=\beta_{n,\lambda}(0)$ are called the degenerate Bernoulli numbers. \\
For $Y=1$, $\beta_{n,\lambda}^{(r)}(x)= \beta_{n,\lambda}^{1,(r)}(x)$ are the degenerate Bernoulli polynomials of order $r$, given by
\begin{equation}
\bigg(\frac{t}{e_{\lambda}(t)-1}\bigg)^{r}e_{\lambda}^{x}(t)=\sum_{n=0}^{\infty}\beta_{n,\lambda}^{(r)}(x)\frac{t^n}{n!}.\label{34a}
\end{equation}
We remark that $\beta_{n,\lambda}(x) \rightarrow B_{n}(x)$, and $\beta_{n,\lambda}^{(r)}(x) \rightarrow B_{n}^{(r)}(x)$, as $\lambda$ tends to $0$. \par
The Frobenius-Euler numbers of order $r$, $H_{n}^{(r)}(u)$, are defined by
\begin{equation}
\bigg(\frac{1-u}{e^{t}-u} \bigg)^{r}=\sum_{n=0}^{\infty}H_{n}^{(r)}(u)\frac{t^{n}}{n!}, \quad(u \ne 1). \label{35a}
\end{equation}
The degenerate Frobenius-Euler numbers of order $r$, $h_{n,\lambda}^{(r)}(u)$, are given by
\begin{equation}
\Big(\frac{1-u}{e_{\lambda}(t)-u} \Big)^{r}=\sum_{n=0}^{\infty} h_{n,\lambda}^{(r)}(u)\frac{t^{n}}{n!}, \quad (u \ne 1). \label{36a}
\end{equation}

\section{Review of umbral calculus}

Here we briefly go over very basic facts about umbral calculus. For more details on this, we recommend the reader to refer to [8,27,28,31]. \par
Let $\mathbb{C}$ be the field of complex numbers. Then $\mathcal{F}$ denotes the algebra of formal power series in $t$ over $\mathbb{C}$, given by
\begin{displaymath}
 \mathcal{F}=\bigg\{f(t)=\sum_{k=0}^{\infty}a_{k}\frac{t^{k}}{k!}~\bigg|~a_{k}\in\mathbb{C}\bigg\},
\end{displaymath}
and $\mathbb{P}=\mathbb{C}[x]$ indicates the algebra of polynomials in $x$ with coefficients in $\mathbb{C}$. \par
Let $\mathbb{P}^{*}$ be the $\mathbb{C}$-vector space of all $\mathbb{C}$- linear functionals on $\mathbb{P}$. If $\langle L|p(x)\rangle$ denotes the action of the $\mathbb{C}$-linear functional $L$ on the polynomial $p(x)$, then the vector space operations on $\mathbb{P}^{*}$ are defined by
\begin{displaymath}
\langle L+M|p(x)\rangle=\langle L|p(x)\rangle+\langle M|p(x)\rangle,\quad\langle cL|p(x)\rangle=c\langle L|p(x)\rangle,
\end{displaymath}
where $c$ is a complex number. \\
For $f(t)\in\mathcal{F}$ with $\displaystyle f(t)=\sum_{k=0}^{\infty}a_{k}\frac{t^{k}}{k!}\displaystyle$, we define the $\mathbb{C}$-linear functional on $\mathbb{P}$ by
\begin{equation}\label{1b}
\langle f(t)|x^{k}\rangle=a_{k}.
\end{equation}
From \eqref{1b}, we note that
\begin{equation*}
 \langle t^{k}|x^{n}\rangle=n!\delta_{n,k},\quad(n,k\ge 0),
\end{equation*}
where $\delta_{n,k}$ is the Kronecker's symbol. \\
Some remarkable $\mathbb{C}$-linear functionals are as follows:
\begin{align}
&\langle e^{yt}|p(x) \rangle=p(y), \nonumber \\
&\langle e^{yt}-1|p(x) \rangle=p(y)-p(0), \label{2b} \\
& \bigg\langle \frac{e^{yt}-1}{t}\bigg |p(x) \bigg\rangle = \int_{0}^{y}p(u) du.\nonumber
\end{align} \par
Let
\begin{equation}\label{3b}
 f_{L}(t)=\sum_{k=0}^{\infty}\langle L|x^{k}\rangle\frac{t^{k}}{k!}.
\end{equation}
Then, by \eqref{1b} and \eqref{3b}, we get
\begin{displaymath}
    \langle f_{L}(t)|x^{n}\rangle=\langle L|x^{n}\rangle.
\end{displaymath}
That is, $f_{L}(t)=L$. Additionally, the map $L\longmapsto f_{L}(t)$ is a $\mathbb{C}$-vector space isomorphism from $\mathbb{P}^{*}$ onto $\mathcal{F}$.\\
Transporting the multiplication in $\mathcal{F}$ to $\mathbb{P}^{*}$ via this isomorphism
gives an algebra structure on $\mathbb{P}^{*}$. This means that the product of $L,\,M \in \mathbb{P}^{*}$ is given by
\begin{equation*}
\langle LM|x^{n} \rangle=\sum_{k=0}^{n}\binom{n}{k}\langle L|x^{k} \rangle \langle M|x^{n-k} \rangle,
\end{equation*}
and the map $L\longmapsto f_{L}(t)$ is now an $\mathbb{C}$-algebra isomorphism.
$\mathcal{F}$ is called umbral algebra which is the algebra of $\mathbb{C}$-linear functionals on $\mathbb{P}$. The umbral calculus is characterized as the study of the umbral algebra. \par
For each nonnegative integer $k$, the differential operator $t^k$ on $\mathbb{P}$ is defined by
\begin{equation}\label{4b}
t^{k}x^n=\left\{\begin{array}{cc}
(n)_{k}x^{n-k}, & \textrm{if $k\le n$,}\\
0, & \textrm{if $k>n$,}
\end{array}\right.
\end{equation}
where $(n)_{k}=n(n-1)\cdots(n-k+1), \,\,(k \ge 1)$, and $(n)_{0}=1$.\\
Extending \eqref{4b} linearly, any power series
\begin{displaymath}
 f(t)=\sum_{k=0}^{\infty}\frac{a_{k}}{k!}t^{k}\in\mathcal{F}
\end{displaymath}
gives the differential operator on $\mathbb{P}$ defined by
\begin{equation}\label{5b}
 f(t)x^n=\sum_{k=0}^{n}\binom{n}{k}a_{k}x^{n-k},\quad(n\ge 0).
\end{equation}
It should be observed that, for any formal power series $f(t)$ and any polynomial $p(x)$, we have
\begin{equation}\label{6b}
\langle f(t) | p(x) \rangle =\langle 1 | f(t)p(x) \rangle =f(t)p(x)|_{x=0}.
\end{equation}
Here we note that an element $f(t)$ of $\mathcal{F}$ is a formal power series, a $\mathbb{C}$-linear functional and a differential  operator. Some notable differential operators are as follows:
\begin{align}
&e^{yt}p(x)=p(x+y), \nonumber\\
&(e^{yt}-1)p(x)=p(x+y)-p(x), \label{7b}\\
&\frac{e^{yt}-1}{t}p(x)=\int_{x}^{x+y}p(u) du.\nonumber
\end{align}

The order $o(f(t))$ of the power series $f(t)(\ne 0)$ is the smallest integer for which $a_{k}$ does not vanish. If $o(f(t))=0$, then $f(t)$ is called an invertible series. If $o(f(t))=1$, then $f(t)$ is called a delta series. \par
For $f(t),g(t)\in\mathcal{F}$ with $o(f(t))=1$ and $o(g(t))=0$, there exists a unique sequence $s_{n}(x)$ of polynomials with $\mathrm{deg}\, s_{n}(x)=n$ such that
\begin{equation} \label{8b}
\big\langle g(t)f(t)^{k}|s_{n}(x)\big\rangle=n!\delta_{n,k},\quad(n,k\ge 0).
\end{equation}
The sequence $s_{n}(x)$ is said to be the Sheffer sequence for $(g(t),f(t))$, which is denoted by $s_{n}(x)\sim (g(t),f(t))$. We observe from \eqref{8b} that
\begin{equation}\label{9b}
s_{n}(x)=\frac{1}{g(t)}p_{n}(x),
\end{equation}
where $p_{n}(x)=g(t)s_{n}(x) \sim (1,f(t))$.\\
In particular, if $s_{n}(x) \sim (g(t),t)$, then $p_{n}(x)=x^n$, and hence
\begin{equation}\label{10b}
s_{n}(x)=\frac{1}{g(t)}x^n.
\end{equation}
It is well known that $s_{n}(x)\sim (g(t),f(t))$ if and only if
\begin{equation}\label{11b}
\frac{1}{g\big(\overline{f}(t)\big)}e^{x\overline{f}(t)}=\sum_{k=0}^{\infty}s_{k}(x)\frac{t^{k}}{k!},
\end{equation}
for all $x\in\mathbb{C}$, where $\overline{f}(t)$ is the compositional inverse of $f(t)$ such that $\overline{f}(f(t))=f(\overline{f}(t))=t$. \par
The following equations \eqref{12b}, \eqref{13b}, and \eqref{14b} are equivalent to the fact that  $s_{n}\left(x\right)$ is Sheffer for $\left(g\left(t\right),f\left(t\right)\right)$, for some invertible $g(t)$:
\begin{align}
f\left(t\right)s_{n}\left(x\right)&=ns_{n-1}\left(x\right),\quad\left(n\ge0\right),\label{12b}\\
s_{n}\left(x+y\right)&=\sum_{j=0}^{n}\binom{n}{j}s_{j}\left(x\right)p_{n-j}\left(y\right),\label{13b}
\end{align}
with $p_{n}\left(x\right)=g\left(t\right)s_{n}\left(x\right),$
\begin{equation}\label{14b}
s_{n}\left(x\right)=\sum_{j=0}^{n}\frac{1}{j!}\big\langle{g\left(\overline{f}\left(t\right)\right)^{-1}\overline{f}\left(t\right)^{j}}\big |{x^{n}\big\rangle}x^{j}.
\end{equation} \par
Let $p_{n}(x),\,q_{n}(x)=\sum_{k=0}^{n}q_{n,k}x^{k}$ be sequences of polynomials. Then the umbral composition of $q_{n}(x)$ with $p_{n}(x)$ is defined to be the sequence
\begin{equation}\label{15b}
q_{n}({\bf{p}}(x))=\sum_{k=0}^{n}q_{n,k}p_{k}(x).
\end{equation}

\section{Representations by probabilistic Bernoulli and degenerate Bernoulli polynomials}
Our interest here is to derive formulas expressing arbitrary polynomials in terms of probabilisitc Bernoulli polynomials associated with $Y$, $B_{n}^{Y}(x)$, and probabilistic degenerate Bernoulli polynomials associated with $Y$, $\beta_{n,\lambda}^{Y}(x)$.  \par
(a) First, we treat the problem of representing arbitrary polynomials by the probabilisitc Bernoulli polynomials associated with $Y$. From \eqref{23a} and \eqref{11b}, we observe that
\begin{align}
&B_{n}^{Y}(x) \sim \bigg( g(t)=\frac{e^{t}-1}{f(t)}, f(t) \bigg), \label{1c} \\
&(x)_{n} \sim (1,\,\, e^{t}-1), \label{2c}
\end{align}
where the compositional inverse of $f(t)$ is given by $\bar{f}(t)=\log E[e^{Yt}]$.
Thus, by \eqref{12b}, we have
\begin{equation}
f(t)B_{n}^{Y}(x)=nB_{n-1}^{Y}(x),\,\,(e^{t}-1)(x)_{n}=n(x)_{n-1}, \label{3c}
\end{equation}
and hence, by \eqref{7b}, we get
\begin{equation}
\Delta(x)_{n}=(e^{t}-1)(x)_{n}=n(x)_{n-1}. \label{4c}
\end{equation}
Here we need to observe that $\bar{f}(t)$ is a delta series. Indeed, one shows that
\begin{equation*}
\bar{f}(t)=E[Y]t+\sum_{n=2}^{\infty}\sum_{j=1}^{n}(-1)^{j-1}(j-1)!S_{2}^{Y}(n,j)\frac{t^{n}}{n!}.
\end{equation*} \par
We note from \eqref{11a}, \eqref{15a} and \eqref{23a} that
\begin{align}
\sum_{n=0}^{\infty}&\big(B_{n}^{Y}(x+1)-B_{n}^{Y}(x) \big)\frac{t^{n}}{n!}=t e^{x \log E[e^{Yt}]} \label{5c}\\
&=t \sum_{j=0}^{\infty}x^{j}\frac{1}{j!}\Big(\log\big(1+(E[e^{Yt}]-1)\big)\Big)^{j} \nonumber \\
&=t \sum_{j=0}^{\infty} x^{j}\sum_{k=j}^{\infty}S_{1}(k,j)\frac{1}{k!}(E[e^{Yt}]-1)^{k}\nonumber\\
&=t \sum_{j=0}^{\infty} x^{j}\sum_{k=j}^{\infty}S_{1}(k,j)\sum_{n=k}^{\infty}S_{2}^{Y}(n,k)\frac{t^{n}}{n!} \nonumber \\
&=\sum_{n=0}^{\infty}\sum_{j=0}^{n}\sum_{k=j}^{n}(n+1)S_{1}(k,j)S_{2}^{Y}(n,k)x^{j}\frac{t^{n+1}}{(n+1)!} \nonumber \\
&=\sum_{n=1}^{\infty}\sum_{j=0}^{n-1}\sum_{k=j}^{n-1}n S_{1}(k,j)S_{2}^{Y}(n-1,k)x^{j}\frac{t^{n}}{n!}. \nonumber
\end{align}
Thus, from \eqref{5c} and \eqref{11a}, we obtain
\begin{align}
B_{n}^{Y}(x+1)-B_{n}^{Y}(x) &=n\sum_{k=0}^{n-1}S_{2}^{Y}(n-1,k)\sum_{j=0}^{k}S_{1}(k,j) x^{j} \label{6c} \\
&=n\sum_{k=0}^{n-1}S_{2}^{Y}(n-1,k) (x)_{k}, \quad (n \ge 1). \nonumber
\end{align}
As $S_{2}^{Y}(n,0)=\delta_{n,0}$ (see \eqref{15a}), from \eqref{6c} we get
\begin{equation}
B_{n}^{Y}(1)-B_{n}^{Y}(0)=nS_{2}^{Y}(n-1,0)=\delta_{n,1}. \label{7c}
\end{equation}
Here $\delta_{n,0}$ and $\delta_{n,1}$ are the Kronecker's delta. \par
Let $p(x) \in \mathbb{C}[x]$ be a polynomial of degree $n$,  $(n \ge 1)$, and let
\begin{equation}
p(x)=\sum_{k=0}^{n}a_{k}B_{k}^{Y}(x). \label{8c}
\end{equation}
Now, we consider
\begin{equation}
a(x)=p(x+1)-p(x)=\Delta p(x)=(e^{t}-1)p(x). \label{9c}
\end{equation}
Then, by \eqref{6c} and \eqref{8c}, we have
\begin{align}
a(x)&=\sum_{k=1}^{n}a_{k}\big(B_{k}^{Y}(x+1)-B_{k}^{Y}(x) \big) \label{10c}\\
&=\sum_{k=1}^{n}k a_{k}\sum_{j=0}^{k-1}S_{2}^{Y}(k-1,j)(x)_{j} \nonumber \\
&=\sum_{k=0}^{n-1}(k+1) a_{k+1}\sum_{j=0}^{k}S_{2}^{Y}(k,j)(x)_{j}, \nonumber
\end{align}
where we used the fact that $B_{0}^{Y}(x)=\frac{1}{E[Y]}$ is a constant polynomial (see \eqref{23a}). \\
By using \eqref{4c} and \eqref{10c}, we have
\begin{align}
\Delta^{r} a(x)&=\sum_{k=r}^{n-1}(k+1) a_{k+1}\sum_{j=r}^{k}S_{2}^{Y}(k,j)\Delta^{r}(x)_{j} \label{11c}\\
&=\sum_{k=r}^{n-1}(k+1) a_{k+1}\sum_{j=r}^{k}S_{2}^{Y}(k,j)(j)_{r}(x)_{j-r}. \nonumber
\end{align}
From \eqref{11c}, we get
\begin{equation}
\frac{1}{r!}\Delta^{r} a(x)|_{x=0}=\sum_{k=r}^{n-1}S_{2}^{Y}(k,r)(k+1)a_{k+1}.\label{12c}
\end{equation}
By inversion (see Proposition 1.1) and from \eqref{12c}, we obtain
\begin{align}
(r+1)a_{r+1}&=\sum_{j=r}^{n-1}S_{1}^{Y}(j,r)\frac{1}{j!} \Delta^{j}a(x)|_{x=0} \label{13c}\\
&=\sum_{j=r}^{n-1}S_{1}^{Y}(j,r)\frac{1}{j!} \Delta^{j+1}p(0). \nonumber
\end{align}
Thus, from \eqref{13c}, we obtain
\begin{equation}
a_{r+1}=\frac{1}{r+1}\sum_{j=r}^{n-1}S_{1}^{Y}(j,r)\frac{1}{j!} \Delta^{j+1}p(0), \quad(r=0,1,\dots, n-1). \label{14c}
\end{equation} \par
We have two alternative expressions for $a_{r+1}$. Firstly, by using \eqref{10a}, \eqref{4c} and \eqref{13c}, we have
\begin{align}
a_{r+1}&=\frac{1}{r+1}\sum_{j=r}^{n-1}S_{1}^{Y}(j,r)\frac{1}{j!} \Delta^{j}a(x)|_{x=0},  \label{15c} \\
&=\frac{1}{r+1}\sum_{j=r}^{n-1}S_{1}^{Y}(j,r)\frac{1}{j!}(e^{t}-1)^{j}a(x)|_{x=0} \nonumber \\
&=\frac{1}{r+1}\sum_{j=r}^{n-1}S_{1}^{Y}(j,r)\sum_{k=j}^{n-1}S_{2}(k,j)\frac{1}{k!}a^{(k)}(x)|_{x=0} \nonumber \\
&=\frac{1}{r+1}\sum_{k=r}^{n-1}\sum_{j=r}^{k}\frac{1}{k!}S_{2}(k,j)S_{1}^{Y}(j,r)\Delta p^{(k)}(0),\quad(r=0,1,\dots n-1), \nonumber
\end{align}
where $p^{(k)}(x)=\big(\frac{d}{dx} \big)^{k}p(x)$. \\
To prceed further, we note from \eqref{15c} that
\begin{equation}
a_{r+1}=\frac{1}{r+1}\sum_{j=r}^{n-1}S_{1}^{Y}(j,r)\frac{1}{j!}\Delta^{j+1}p(x)|_{x=0}. \label{16c}
\end{equation}
Secondly, from \eqref{30a} and \eqref{16c}, we have
\begin{equation}
a_{r+1}=\frac{1}{r+1}\sum_{j=r}^{n-1}\sum_{k=0}^{j+1}(-1)^{j+1-k}\frac{1}{j!}\binom{j+1}{k}S_{1}^{Y}(j,r)p(k), \quad(r=0,1,\dots n-1).\label{17c}
\end{equation} \par
It still remains to determine $a_{0}$. To do this, from \eqref{2b}, \eqref{1c}, \eqref{3c} and \eqref{7c}, we observe first that
\begin{align}
\langle g(t)|p(x) \rangle&=\Big\langle \frac{e^{t}-1}{f(t)}\Big \vert p(x)  \Big\rangle \label{18c} \\
&=\sum_{k=0}^{n} a_{k}\Big\langle \frac{e^{t}-1}{f(t)}\Big \vert B_{k}^{Y}(x) \Big\rangle\nonumber \\
&=\sum_{k=0}^{n} a_{k}\Big\langle \frac{e^{t}-1}{f(t)}\Big \vert f(t) \frac{1}{k+1}B_{k+1}^{Y}(x) \Big\rangle\nonumber \\
&=\sum_{k=0}^{n} \frac{a_{k}}{k+1}\langle e^{t}-1|B_{k+1}^{Y}(x) \rangle \nonumber \\
&=\sum_{k=0}^{n} \frac{a_{k}}{k+1}\big(B_{k+1}^{Y}(1)-B_{k+1}^{Y}(0) \big) \nonumber\\
&=\sum_{k=0}^{n}\frac{a_{k}}{k+1} \delta_{k+1,1}=a_{0}. \nonumber
\end{align}
As $g(t)=\frac{e^{t}-1}{t} \frac{t}{f(t)}$, from \eqref{2b} and \eqref{18c}, we obtain
\begin{equation}
a_{0}=\Big \langle \frac{e^{t}-1}{t} \Big \vert \frac{t}{f(t)} p(x) \Big \rangle=\int_{0}^{1} \frac{t}{f(t)}p(x) dx. \label{19c}
\end{equation}
Now, we obtain the following theorem from \eqref{14c}, \eqref{15c}, \eqref{17c} and \eqref{19c}.
\begin{theorem}
Let $p(x) \in \mathbb{C}[x]$ be a polynomial of degree $n$, and let $p(x)=\sum_{k=0}^{n}a_{k}B_{k}^{Y}(x)$. Then we have
\end{theorem}
\begin{align*}
&a_{0}=\int_{0}^{1} \frac{t}{f(t)}p(x) dx, \\
&a_{r+1}=\frac{1}{r+1}\sum_{j=r}^{n-1}S_{1}^{Y}(j,r)\frac{1}{j!} \Delta^{j+1}p(0) \\
&=\frac{1}{r+1}\sum_{k=r}^{n-1}\sum_{j=r}^{k}\frac{1}{k!}S_{2}(k,j)S_{1}^{Y}(j,r)\Delta p^{(k)}(0) \\
&=\frac{1}{r+1}\sum_{j=r}^{n-1}\sum_{k=0}^{j+1}(-1)^{j+1-k}\frac{1}{j!}\binom{j+1}{k}S_{1}^{Y}(j,r)p(k), \quad(r=0,1,\dots n-1),
\end{align*}
where $f(t)$ is the compositional inverse of $\bar{f}(t)=\log E[e^{Yt}]$.
\begin{remark}
Consider the case $Y=1$. From Theorem 3.1, for $r=0,1,\dots,n-1$, we have
\end{remark}
\begin{align}
a_{r+1}&=\frac{1}{r+1}\sum_{k=r}^{n-1}\frac{1}{k!}(p^{(k)}(1)-p^{(k)}(0))\sum_{j=r}^{k}S_{2}(k,j)S_{1}(j,r) \label{20c}\\
&=\frac{1}{r+1}\sum_{k=r}^{n-1}\frac{1}{k!}(p^{(k)}(1)-p^{(k)}(0))\delta_{k,r}\nonumber \\
&=\frac{1}{(r+1)!}(p^{(r)}(1)-p^{(r)}(0)) \nonumber \\
&=\frac{1}{(r+1)!}\int_{0}^{1}p^{(r+1)}(x) dx. \nonumber
\end{align}
As $f(t)=t$, from Theorem 3.1, we get
\begin{equation}
a_{0}=\int_{0}^{1}p(x) dx. \label{21c}
\end{equation}
So, by \eqref{20c} and \eqref{21c}, we recover the well-known formula:
\begin{equation}
p(x)=\sum_{k=0}^{n}a_{k}B_{k}(x),\quad a_{k}=\frac{1}{k!}\int_{0}^{1}p^{(k)}(x) dx,\quad(k=0,1,\dots,n). \label{21-1c}
\end{equation} \par
(b) Next, we consider the problem of representing arbitrary polynomials by the probabilisitc degenerate Bernoulli polynomials associated with $Y$. From \eqref{25a} and \eqref{11b}, we observe that
\begin{equation}
\beta_{n,\lambda}^{Y}(x) \sim \bigg( g(t)=\frac{e^{t}-1}{f(t)}, f(t) \bigg),\label{22c}
\end{equation}
where the compositional inverse of $f(t)$ is given by $\bar{f}(t)=\log E[e_{\lambda}^{Y}(t)]$.
Thus, by \eqref{11b}, we have
\begin{equation}
f(t)\beta_{n,\lambda}^{Y}(x)=n\beta_{n-1,\lambda}^{Y}(x). \label{23c}
\end{equation}
Here we have to see that $\bar{f}(t)$ is a delta series. Indeed, we have
\begin{equation}
\bar{f}(t)=E[Y]t+\sum_{n=2}^{\infty}\sum_{j=1}^{n}(-1)^{j-1}(j-1)!S_{2,\lambda}^{Y}(n,j) \frac{t^{n}}{n!}. \label{24c}
\end{equation}
By using the definition in \eqref{25a} and proceeding just as in \eqref{6c}, we can show that
\begin{equation}
\beta_{n,\lambda}^{Y}(x+1)-\beta_{n,\lambda}^{Y}(x)=n \sum_{k=0}^{n-1}S_{2,\lambda}^{Y}(n-1,k)(x)_{k}, \,\, (n \ge 1). \label{25c}
\end{equation}
From \eqref{25c} and as $S_{2}^{Y}(n,0)=\delta_{n,0}$, it follows that
\begin{equation}
\beta_{n,\lambda}^{Y}(1)-\beta_{n,\lambda}^{Y}(0)=\delta_{n,1}. \label{26c}
\end{equation} \par
Let $p(x) \in \mathbb{C}[x]$ be a polynomial of degree $n$,  $(n \ge 1)$, and let
\begin{equation}
p(x)=\sum_{k=0}^{n}a_{k}\beta_{k,\lambda}^{Y}(x). \label{27c}
\end{equation}
Now, we consider
\begin{equation}
a(x)=p(x+1)-p(x)=\Delta p(x)=(e^{t}-1)p(x). \label{28c}
\end{equation}
Then, by proceeding just as in \eqref{10c}, \eqref{11c} and \eqref{12c}, we have
\begin{equation}
\frac{1}{r!}\Delta^{r} a(x)|_{x=0}=\sum_{k=r}^{n-1}S_{2,\lambda}^{Y}(k,r)(k+1)a_{k+1},\,\,\,(r=0,1,\dots,n-1). \label{29c}
\end{equation}
By inversion (see Proposition 1.2) and from \eqref{29c}, we obtain
\begin{align}
a_{r+1}&=\frac{1}{r+1}\sum_{j=r}^{n-1}S_{1,\lambda}^{Y}(j,r)\frac{1}{j!}\Delta^{j}a(x)|_{x=0} \label {30c}\\
&=\frac{1}{r+1}\sum_{j=r}^{n-1}S_{1,\lambda}^{Y}(j,r)\frac{1}{j!}\Delta^{j+1}p(0). \nonumber
\end{align} \par
We have two alternative expressions for $a_{r+1}$, which are given by
\begin{align}
a_{r+1}&=\frac{1}{r+1}\sum_{k=r}^{n-1}\sum_{j=r}^{k}\frac{1}{k!}S_{2}(k,j)S_{1,\lambda}^{Y}(j,r)\Delta p^{(k)}(0) \label{31c} \\
&=\frac{1}{r+1}\sum_{j=r}^{n-1}\sum_{k=0}^{j+1}(-1)^{j+1-k}\frac{1}{j!}\binom{j+1}{k}S_{1,\lambda}^{Y}(j,r)p(k), \quad(r=0,1,\dots n-1), \nonumber
\end{align}
where $p^{(k)}(x)=\big(\frac{d}{dx} \big)^{k}p(x)$.
In addition, by proceeding just as in \eqref{18c} and \eqref{19c}, we have
\begin{equation}
a_{0}=\int_{0}^{1}\frac{t}{f(t)}p(x) dx. \label{32c}
\end{equation}
\begin{theorem}
Let $p(x) \in \mathbb{C}[x]$ be a polynomial of degree $n$, and let $p(x)=\sum_{k=0}^{n}a_{k}\beta_{k,\lambda}^{Y}(x)$. Then we have
\end{theorem}
\begin{align*}
&a_{0}=\int_{0}^{1} \frac{t}{f(t)}p(x) dx, \\
&a_{r+1}=\frac{1}{r+1}\sum_{j=r}^{n-1}S_{1,\lambda}^{Y}(j,r)\frac{1}{j!} \Delta^{j+1}p(0) \\
&=\frac{1}{r+1}\sum_{k=r}^{n-1}\sum_{j=r}^{k}\frac{1}{k!}S_{2}(k,j)S_{1,\lambda}^{Y}(j,r)\Delta p^{(k)}(0) \\
&=\frac{1}{r+1}\sum_{j=r}^{n-1}\sum_{k=0}^{j+1}(-1)^{j+1-k}\frac{1}{j!}\binom{j+1}{k}S_{1,\lambda}^{Y}(j,r)p(k), \quad(r=0,1,\dots n-1),
\end{align*}
where $f(t)$ is the compositional inverse of $\bar{f}(t)=\log E[e_{\lambda}^{Y}(t)]$. \par
As $\frac{\lambda t}{e^{\lambda t}-1}e^{x t}=\sum_{n=0}^{\infty}\lambda^{n}B_{n}\big(\frac{x}{\lambda}\big) \frac{t^{n}}{n!}$, we have
\begin{equation}
\lambda^{n}B_{n}\Big(\frac{x}{\lambda} \Big)=\frac{\lambda t}{e^{\lambda t}-1} x^{n}. \label{33c}
\end{equation}
From Theorem 3.3 and \eqref{33c}, we get the following theorem.
\begin{remark}
Let $Y=1$. Then $f(t)=\frac{1}{\lambda}(e^{\lambda t}-1)$. Now, $p(x)=\sum_{k=0}^{n}a_{k}\beta_{k,\lambda}(x)$, where
\begin{align}
&a_{0}=\int_{0}^{1} \frac{\lambda t}{e^{\lambda t}-1}p(x) dx=\int_{0}^{1} p\Big(\lambda \mathbf{B} \Big(\frac{x}{\lambda}\Big)\Big) dx, \label{34c}\\
&a_{r+1}=\frac{1}{r+1}\sum_{j=r}^{n-1}S_{1,\lambda}(j,r)\frac{1}{j!} \Delta^{j+1}p(0) \nonumber \\
&=\frac{1}{r+1}\sum_{k=r}^{n-1}\sum_{j=r}^{k}\frac{1}{k!}S_{2}(k,j)S_{1,\lambda}(j,r)\Delta p^{(k)}(0) \nonumber \\
&=\frac{1}{r+1}\sum_{j=r}^{n-1}\sum_{k=0}^{j+1}(-1)^{j+1-k}\frac{1}{j!}\binom{j+1}{k}S_{1,\lambda}(j,r)p(k), \quad(r=0,1,\dots n-1). \nonumber
\end{align}
Here $p\big(\lambda \mathbf{B} \big(\frac{x}{\lambda}\big)\big)$ is the umbral composition of $p(x)$ with $\lambda^{i}B_{i}\big(\frac{x}{\lambda}\big)$. \par
In [16], the following different expressions for $a_{r+1},\,(r=0,1,\dots n-1)$,\, are derived.
\begin{align}
a_{r+1}&=\frac{1}{(r+1)! \lambda^{r}}\Delta_{\lambda}^{r} \Delta p(0) \label{35c} \\
&=\frac{1}{(r+1)! \lambda^{r}}\sum_{j=0}^{r}\binom{r}{j}(-1)^{r-j} \Delta p(j \lambda) \nonumber \\
&=\frac{1}{r+1}\sum_{l=r}^{n}S_{2}(l,r)\frac{\lambda^{l-r}}{l!} \Delta p^{(l)}(0). \nonumber
\end{align}
\end{remark}

\section{Representations by higher-order probabilistic Bernoulli and degenerate Bernoulli polynomials}

Our interest here is to derive formulas expressing arbitrary polynomials in terms of higher-order probabilistic Bernoulli polynomials associated with $Y$, $B_{n}^{Y,(r)}(x)$, and higher-order probabilistic degenerate Bernoulli polynomials associated with $Y$, $\beta_{n, \lambda}^{Y,(r)}(x)$. \par
(a) First, we treat the problem of representing arbitrary polynomials by the higher-order  probabilisitc Bernoulli polynomials associated with $Y$. \par
We note from \eqref{32a} and \eqref{11b} that
\begin{equation}
B_{n}^{Y,(r)}(x) \sim (g(t)^{r}, f(t)), \label{1d}
\end{equation}
where $g(t)=\frac{e^{t}-1}{f(t)}$, \,and the compositional inverse of $f(t)$ is given by
$\bar{f}(t)=\log E[e^{Yt}]$. \\
So, by \eqref{12b}, we have
\begin{equation}
f(t)B_{n}^{Y,(r)}(x)=n B_{n-1}^{Y,(r)}(x). \label{2d}
\end{equation}
We observe from \eqref{32a} that
\begin{equation}
\Delta B_{n}^{Y,(r)}(x)=B_{n}^{Y,(r)}(x+1)-B_{n}^{Y, (r)}(x)=nB_{n-1}^{Y,(r-1)}(x),\,\,(n \ge 0). \label{3d}
\end{equation}
We see from \eqref{7b}, \eqref {1d}, \eqref{2d} and \eqref{3d} that
\begin{align}
g(t)B_{n}^{Y,(r)}(x)&=\frac{e^{t}-1}{f(t)} f(t)\frac{B_{n+1}^{Y,(r)}(x)}{n+1}=\frac{1}{n+1}(e^{t}-1)B_{n+1}^{Y,(r)}(x) \label{4d} \\
&=\frac{1}{n+1}\big(B_{n+1}^{Y,(r)}(x+1)-B_{n+1}^{Y,(r)}(x) \big)=B_{n}^{Y,(r-1)}(x).\nonumber
\end{align} \par
Let $p(x) \in \mathbb{C}[x]$ be a polynomial of degree $n$, and write
\begin{equation}
p(x)=\sum_{k=0}^{n}a_{k}B_{k}^{Y, (r)}(x). \label{5d}
\end{equation}
As $ \langle g(t)^{r}f(t)^{k} | B_{n}^{Y,(r)}(x) \rangle =n! \delta_{n,k}$ (see \eqref{8b}), we have
\begin{equation}
a_{k}=\frac{1}{k!} \langle f(t)^{k}g(t)^{r} | p(x) \rangle. \label{6d}
\end{equation}
Let $a$ be any nonnegative integer. Then we have
\begin{equation}
f(t)^{a}p(x)=\Big(\frac{f(t)}{t}\Big)^{a} t^{a}p(x)=\Big(\frac{f(t)}{t}\Big)^{a}p^{(a)}(x). \label{7d}
\end{equation}
Let $I$ be the linear operator defined on $\mathbb{P}$, which is given by (see \eqref{7b})
\begin{equation}
Iq(x)=\frac{e^{t}-1}{t}q(x)=\int_{x}^{x+1}q(u) du. \label{8d}
\end{equation}
Then, noting that $g(t)=\frac{e^{t}-1}{t} \frac{t}{f(t)}$, we obtain the following expression:
\begin{equation}
g(t)^{a}p(x)=I^{a}\Big(\frac{t}{f(t)}\Big)^{a}p(x). \label{9d}
\end{equation}
Using $\big(\frac{e^{t}-1}{t}\big)^{m}=\sum_{l=0}^{\infty}S_{2}(l+m,m)\frac{m!}{(l+m)!}t^{l}$ and \eqref{9d}, we get another expression:
\begin{equation}
g(t)^{a}p(x)=\sum_{l=0}^{n}S_{2}(l+a,a)\frac{(a)!}{(l+a)!}\Big(\frac{t}{f(t)}\Big)^{a}p^{(l)}(x). \label{10d}
\end{equation}
Noting $f(t)g(t)=e^{t}-1$, cosidering \eqref{6d} for the cases $r>n$ and $r \le n$ separately (see [16]) and from \eqref{7d}, \eqref{9d} and \eqref{10d}, we obtain the following theorem.
\begin{theorem}
Let $p(x) \in \mathbb{C}[x]$ be a polynomial of degree $n$, and let $g(t)=\frac{e^{t}-1}{f(t)}$, with $\bar{f}(t)=\log E[e^{Y t}]$. Then we have the following: \\
(a) For $r > n$, we have
\begin{equation*}
p(x)=\sum_{k=0}^{n}\frac{1}{k!}\sum_{j=0}^{k}(-1)^{k-j}\binom{k}{j}\big(g(t)^{r-k}p(j)\big)B_{k}^{Y,(r)}(x).
\end{equation*}
(b) For $r \le n$, we have
\begin{align*}
p(x)=\sum_{k=0}^{r-1}\frac{1}{k!}\sum_{j=0}^{k}(-1)^{k-j}\binom{k}{j}\big(g(t)^{r-k}p(j)\big)B_{k}^{Y,(r)}(x) \\
\quad \quad +\sum_{k=r}^{n}\frac{1}{k!}\sum_{j=0}^{r}(-1)^{r-j}\binom{r}{j}\big(f(t)^{k-r}p(j)\big)B_{k}^{Y,(r)}(x).
\end{align*}
Here $f(t)^{k-r}p(j)$ and $g(t)^{r-k}p(j)$ can be expressed by
\begin{equation*}
f(t)^{k-r}p(j)=\Big(\frac{f(t)}{t}\Big)^{k-r}p^{(k-r)}(x)\Big\vert_{x=j},
\end{equation*}
\begin{align*}
g(t)^{r-k}p(j)&=I^{r-k}\Big(\frac{t}{f(t)}\Big)^{r-k}p(x)\Big\vert_{x=j} \\
&=\sum_{l=0}^{n}S_{2}(l+r-k,r-k)\frac{(r-k)!}{(l+r-k)!}\Big(\frac{t}{f(t)}\Big)^{r-k}p^{(l)}(x)\Big\vert_{x=j},
\end{align*}
where $I$ is the linear operator on $\mathbb{P}$ given by $Iq(x)=\int_{x}^{x+1}q(u) du$.
\end{theorem}

(b) Next, we treat the problem of representing arbitrary polynomials by the higher-order  probabilisitc degenerate Bernoulli polynomials associated with $Y$. \par
We note from \eqref{34a} and \eqref{11b} that
\begin{equation}
\beta_{n,\lambda}^{Y,(r)}(x) \sim (g(t)^{r}, f(t)), \label{11d}
\end{equation}
where $g(t)=\frac{e^t -1}{f(t)}$, and the compositional inverse of $f(t)$ is given by $\bar{f}(t)=\log E[e_{\lambda}^{Y}(t)]$.
Thus, by \eqref{12b}, we have
\begin{equation}
f(t)\beta_{n,\lambda}^{Y,(r)}(x)=n\beta_{n-1,\lambda}^{Y,(r)}(x)    \label{12d}
\end{equation}
We also note from \eqref{34a} that
\begin{equation}
\Delta \beta_{n,\lambda}^{Y,(r)}(x)=\beta_{n,\lambda}^{Y,(r)}(x+1)-\beta_{n,\lambda}^{Y, (r)}(x)=n\beta_{n-1,\lambda}^{Y,(r-1)}(x),\,\,(n \ge 0). \label{13d}
\end{equation}
Using \eqref{11d}, \eqref{12d} and \eqref{13d}, and proceeding as in \eqref{4d}, we can show that
\begin{equation}
g(t)\beta_{n,\lambda}^{Y,(r)}(x)=\beta_{n,\lambda}^{Y,(r-1)}(x), \label{14d}
\end{equation} \par
Let $p(x) \in \mathbb{C}[x]$ be a polynomial of degree $n$, and write \\
\begin{equation}
p(x)=\sum_{k=0}^{n}a_{k}\beta_{k,\lambda}^{Y, (r)}(x). \label{15d}
\end{equation}
As $ \langle g(t)^{r}f(t)^{k} | \beta_{n,\lambda}^{Y,(r)}(x) \rangle =n! \delta_{n,k}$ (see \eqref{8b}), we have
\begin{equation}
a_{k}=\frac{1}{k!} \langle f(t)^{k}g(t)^{r} | p(x) \rangle. \label{16d}
\end{equation}
Noting $f(t)g(t)=e^{t}-1$, cosidering the cases $r>n$ and $r \le n$ separately and proceeding just as in the case (a), we obtain
\begin{theorem}
Let $p(x) \in \mathbb{C}[x]$ be a polynomial of degree $n$, and let $g(t)=\frac{e^{t}-1}{f(t)}$, with $\bar{f}(t)=\log E[e_{\lambda}^{Y}(t)]$. Then we have the following: \\
(a) For $r > n$, we have
\begin{equation*}
p(x)=\sum_{k=0}^{n}\frac{1}{k!}\sum_{j=0}^{k}(-1)^{k-j}\binom{k}{j}\big(g(t)^{r-k}p(j)\big)\beta_{k,\lambda}^{Y,(r)}(x).
\end{equation*}
(b) For $r \le n$, we have
\begin{align*}
p(x)=\sum_{k=0}^{r-1}\frac{1}{k!}\sum_{j=0}^{k}(-1)^{k-j}\binom{k}{j}\big(g(t)^{r-k}p(j)\big)\beta_{k,\lambda}^{Y,(r)}(x) \\
\quad \quad +\sum_{k=r}^{n}\frac{1}{k!}\sum_{j=0}^{r}(-1)^{r-j}\binom{r}{j}\big(f(t)^{k-r}p(j)\big)\beta_{k,\lambda}^{Y,(r)}(x).
\end{align*}
Here $f(t)^{k-r}p(j)$ and $g(t)^{r-k}p(j)$ can be expressed by
\begin{equation*}
f(t)^{k-r}p(j)=\Big(\frac{f(t)}{t}\Big)^{k-r}p^{(k-r)}(x)\Big\vert_{x=j},
\end{equation*}
\begin{align*}
g(t)^{r-k}p(j)&=I^{r-k}\Big(\frac{t}{f(t)}\Big)^{r-k}p(x)\Big\vert_{x=j} \\
&=\sum_{l=0}^{n}S_{2}(l+r-k,r-k)\frac{(r-k)!}{(l+r-k)!}\Big(\frac{t}{f(t)}\Big)^{r-k}p^{(l)}(x)\Big\vert_{x=j},
\end{align*}
where $I$ is the linear operator on $\mathbb{P}$ given by $Iq(x)=\int_{x}^{x+1}q(u) du$.
\end{theorem}

Assume now that $Y=1$. Then we have
\begin{equation}
f(t)=\frac{1}{\lambda}\big(e^{\lambda t}-1 \big), \quad g(t)=\frac{e^{t}-1}{t}\frac{\lambda t}{e^{\lambda t}-1}. \label{17d}
\end{equation}
As $\big(\frac{\lambda t}{e^{\lambda t}-1} \big)^{a}e^{x t}=\sum_{n=0}^{\infty}\lambda^{n}B_{n}^{(a)}\big(\frac{x}{\lambda}\big) \frac{t^{n}}{n!}$, for any nonnegative integer $a$, we have
\begin{equation}
\lambda^{n}B_{n}^{(a)}\big(\frac{x}{\lambda} \big)=\bigg(\frac{\lambda t}{e^{\lambda t}-1} \bigg)^{a} x^{n}, \label{18d}
\end{equation}
and hence, from \eqref{17d} and \eqref{18d}, we see that
\begin{equation}
\Big(\frac{t}{f(t)}\Big)^{r-k}p(x)=\bigg(\frac{\lambda t}{e^{\lambda t}-1} \bigg)^{r-k}p(x)=p\big(\lambda \mathbf{B}^{(r-k)}(\frac{x}{\lambda})\big), \label{19d}
\end{equation}
where $p\big(\lambda \mathbf{B}^{(r-k)}(\frac{x}{\lambda})\big)$ is the umbral composition of $p(x)$ with $\lambda^{i}B_{i}^{(r-k)}(\frac{x}{\lambda})$.
Thus, from Theorem 4.2 and \eqref{19d}, we have
\begin{align}
g(t)^{r-k}p(x)&=I^{r-k}p\big(\lambda \mathbf{B}^{(r-k)}(\frac{x}{\lambda})\big) \label{20d} \\
&=\sum_{l=0}^{n}S_{2}(l+r-k,r-k)\frac{(r-k)!}{(l+r-k)!}p^{(l)}\big(\lambda \mathbf{B}^{(r-k)}(\frac{x}{\lambda})\big). \nonumber
\end{align}
From \eqref{7b} and \eqref{17d}, we have
\begin{align}
f(t)^{k-r}p(x)&=\Big(\frac{f(t)}{t}\Big)^{k-r}p^{(k-r)}(x) \label{21d} \\
&=\frac{1}{\lambda^{k-r}}\Big(\frac{e^{\lambda t}-1}{t} \Big)^{k-r}p^{(k-r)}(x) \nonumber \\
&=\frac{1}{\lambda^{k-r}}I_{\lambda}^{k-r}p^{(k-r)}(x), \nonumber
\end{align}
where $I_{\lambda}$ is the linear operator on $\mathbb{P}$ defined by
\begin{equation}
I_{\lambda}q(x)=\frac{e^{\lambda t}-1}{t}q(x)=\int_{x}^{x+\lambda}q(u) du. \label{22d}
\end{equation}
Furthermore, we have
\begin{align}
f(t)^{k-r}p(x)&=\frac{1}{\lambda^{k-r}}(e^{\lambda t}-1)^{k-r}p(x) \label{23d} \\
&=\frac{1}{\lambda^{k-r}}\Delta_{\lambda}^{k-r}p(x) \nonumber \\
&=\frac{1}{\lambda^{k-r}}\sum_{l=0}^{k-r}\binom{k-r}{l}(-1)^{k-r-l}p(x+l \lambda), \nonumber
\end{align}
and
\begin{align}
f(t)^{k-r}p(x)&=\frac{(k-r)!}{\lambda^{k-r}}\frac{1}{(k-r)!}(e^{\lambda t}-1)^{k-r}p(x) \label{24d} \\
&=\frac{(k-r)!}{\lambda^{k-r}}\sum_{l=k-r}^{n}S_{2}(l,k-r)\frac{1}{l!}p^{(l)}(x). \nonumber
\end{align}

\begin{remark}
Let $p(x) \in \mathbb{C}[x], with \,\,\mathrm{deg}\, p(x)=n$. Write $p(x)=\sum_{k=0}^{n}a_kB_{k}^{(r)}(x)$. As $\lambda$ tends to $0$, $g(t) \rightarrow \frac{e^t -1}{t},\,\,f(t) \rightarrow t,\,\, p\big(\lambda {\bf{B}}^{(r-k)}\big(\frac{x}{\lambda}\big)\big) \rightarrow p(x)$. Thus, from Theorem 4.2, we recover the following result obtained in \textnormal{[16]}: \par
(a) For $r>n$, we have
\begin{equation*}
p(x)=\sum_{k=0}^{n}\bigg(\sum_{j=0}^{k}\frac{1}{k!}(-1)^{k-j}\binom{k}{j}I^{r-k}p(j)\bigg)B_{k}^{(r)}(x).
\end{equation*}
(b) For $r \leq n$, we have
\begin{align*}
p(x)=&\sum_{k=0}^{r-1}\bigg(\sum_{j=0}^{k}\frac{1}{k!}(-1)^{k-j}\binom{k}{j}I^{r-k}p(j)\bigg)B_{k}^{(r)}(x)\\
&+\sum_{k=r}^{n}\bigg(\sum_{j=0}^{r}\frac{1}{k!}(-1)^{r-j}p^{(k-r)}(j)\bigg)B_{k}^{(r)}(x).
\end{align*}
Here $I$ is the integral operator $Iq(x)=\int_{x}^{x+1}q(u)\,du$.
\end{remark}

\section{Examples}
Here we express $x^{n}$ as linear combinations of probabilistic Bernoulli polynomials associated with $Y$, $B_{k}^{Y}(x)$, and probabilistic degenerate Bernoulli polynomials associated with $Y$, $\beta_{k,\lambda}^{Y}(x)$, for several discrete and continuous random variables $Y$. For those random variables $Y$, we need the explicit computations in [14], for the following:
\begin{equation}
f_{Y}(t),\quad f_{Y,\lambda}(t), \quad S_{1}^{Y}(n,k), \quad S_{1,\lambda}^{Y}(n,k), \label{1e}
\end{equation}
where the compositional inverses of $f_{Y}(t)$ and $f_{Y,\lambda}(t)$ are respectively given by
\begin{equation*}
\bar{f}_{Y}(t)=\log E[e^{Y t}], \quad \bar{f}_{Y,\lambda}(t)=\log E[e_{\lambda}^{Y}(t)].
\end{equation*}
{\it{Note here that we denote the $f(t)$ in Theorem 3.1 and the $f(t)$ in Theorem 3.3 respectively by $f_{Y}(t)$ and $f_{Y,\lambda}(t)$.}}
We start with the following lemma.
\begin{lemma}
For integer $n \ge 0$, we have
\begin{equation*}
\int_{0}^{1}B_{n}(x) dx =\delta_{n,0},\quad \int_{0}^{1} B_{n}(-x) dx=(-1)^{n},
\end{equation*}
where $\delta_{n,0}$ is the Kronecker's delta.
\begin{proof}
The first one follows from the facts $\frac{d}{dx} B_{n}(x)=nB_{n-1}(x)$ and $B_{n}(1)-B_{n}=\delta_{n,1}$, the latter of which follows from \eqref{6c} with $Y=1$. For the second one, we note the following:
\begin{equation}
B_{n+1}(-x)=(-1)^{n+1}\big(B_{n+1}(x)+(n+1)x^{n}\big). \label{2e}
\end{equation}
Now, by using \eqref{2e} and the fact $B_{n+1}(1)-B_{n+1}=\delta_{n,0}$, we see that
\begin{align}
\int_{0}^{1}B_{n}(-x) dx &=\Big[-\frac{1}{n+1}B_{n+1}(-x)\Big]_{0}^{1} \label{3e}\\
&=-\frac{1}{n+1}\Big\{(-1)^{n+1}\delta_{n,0}+(-1)^{n+1}B_{n+1}+(-1)^{n+1}(n+1)-B_{n+1}\Big\}. \nonumber
\end{align}
Now, the result follows from \eqref{3e} by considering the cases that $n$ is odd, $n$ is even with $n \ge 2$, and $n=0$, separately.
\end{proof}
\end{lemma}

\begin{remark}
Shorter proofs of Lemma 5.1 are as follows:
\begin{align*}
&\sum_{n=0}^{\infty}\frac{t^{n}}{n!}\int_{0}^{1}B_{n}(x) dx=\frac{t}{e^{t}-1}\int_{0}^{1}e^{xt} dx=1,\\
&\sum_{n=0}^{\infty}\frac{t^{n}}{n!}\int_{0}^{1}B_{n}(-x) dx=\frac{t}{e^{t}-1}\int_{0}^{1}e^{-xt} dx=e^{-t}=\sum_{n=0}^{\infty}\frac{t^{n}}{n!}(-1)^{n}.
\end{align*}
\end{remark}

(a) Let $Y$ be the Bernoulli random variable. Then the probability mass function of $Y$ is given by (see [29])
\begin{equation*}
p(0)=1-p,\quad p(1)=p, \quad (0< p \le 1).
\end{equation*}
Then, from [14], we have
\begin{align}
&f_{Y}(t)=\log\big(1+\frac{1}{p}(e^{t}-1)\big),\,\,f_{Y,\lambda}(t)=\log_{\lambda}\big(1+\frac{1}{p}(e^{t}-1)\big), \label{4e} \\
&S_{1}^{Y}(n,k)=\frac{1}{p^{n}}S_{1}(n,k),\,\, S_{1,\lambda}^{Y}(n,k)=\frac{1}{p^{n}}S_{1,\lambda}(n,k). \nonumber
\end{align} \par
Let $x^n=\sum_{k=0}^{n}a_{k}B_{k}^{Y}(x)$. We first observe that
\begin{align}
\frac{t}{\log \big(1+\frac{1}{p}(e^{t}-1)\big)}e^{xt}&=\frac{\frac{1}{p}(e^{t}-1)}{\log\big(1+\frac{1}{p}(e^{t}-1)\big)}\frac{t}{\frac{1}{p}(e^{t}-1)}e^{xt} \label{5e} \\
&=\sum_{l=0}^{\infty}b_{l}\frac{1}{p^{l-1}}\frac{1}{l!}(e^{t}-1)^{l}\sum_{m=0}^{\infty}B_{m}(x)\frac{t^{m}}{m!} \nonumber\\
&=\sum_{r=0}^{\infty}\sum_{l=0}^{r}b_{l}\frac{1}{p^{l-1}}S_{2}(r,l)\frac{t^{r}}{r!}\sum_{m=0}^{\infty}B_{m}(x)\frac{t^{m}}{m!} \nonumber \\
&=\sum_{n=0}^{\infty}\sum_{r=0}^{n}\sum_{l=0}^{r}\binom{n}{r}\frac{1}{p^{l-1}}S_{2}(r,l)b_{l}B_{n-r}(x)\frac{t^{n}}{n!}, \nonumber
\end{align}
where $b_{l}$ are the Bernoulli numbers of the second defined by
\begin{equation}
\frac{t}{\log (1+t)}=\sum_{l=0}^{\infty}b_{l} \frac{t^{l}}{l!}. \label{6e}
\end{equation}
Thus, from \eqref{4e} and \eqref{5e}, we get
\begin{equation}
\frac{t}{f_{Y}(t)}x^{n}=\sum_{r=0}^{n}\sum_{l=0}^{r}\binom{n}{r}\frac{1}{p^{l-1}}S_{2}(r,l)b_{l}B_{n-r}(x). \label{7e}
\end{equation}
Now, from Theorem 3.1, \eqref{7e} and Lemma 5.1, we have
\begin{align}
a_{0}&=\int_{0}^{1}\frac{t}{f_{Y}(t)}x^{n} =\sum_{r=0}^{n}\sum_{l=0}^{r}\binom{n}{r}\frac{1}{p^{l-1}}S_{2}(r,l)b_{l}\int_{0}^{1}B_{n-r}(x) dx \label{8e} \\
&=\sum_{r=0}^{n}\sum_{l=0}^{r}\binom{n}{r}\frac{1}{p^{l-1}}S_{2}(r,l)b_{l}\delta_{n,r}=\sum_{l=0}^{n}\frac{1}{p^{l-1}}S_{2}(n,l)b_{l}. \nonumber
\end{align}
From Theorem 3.1, \eqref{4e} and \eqref{8e}, we derive the following representation
\begin{align*}
x^{n}&=\sum_{l=0}^{n}\frac{1}{p^{l-1}}S_{2}(n,l)b_{l}B_{0}^{Y}(x) \\
&\quad\quad +\sum_{k=1}^{n}\Big\{\frac{1}{k}\sum_{j=k-1}^{n-1}\frac{1}{p^{j}}
S_{1}(j,k-1)\frac{1}{j!}\Delta^{j+1}0^{n} \Big\}B_{k}^{Y}(x).
\end{align*} \par
Next, we let $x^{n}= \sum_{k=0}^{n}a_{k}\beta_{k,\lambda}^{Y}(x)$. Proceeding just as before, we get
\begin{equation}
\frac{t}{f_{Y,\lambda}(t)}x^{n}=\sum_{r=0}^{n}\sum_{l=0}^{r}\binom{n}{r}\frac{1}{p^{l-1}}S_{2}(r,l)b_{l,\lambda}B_{n-r}(x), \label{9e}
\end{equation}
where $b_{l,\lambda}$ are the degenerate Bernoulli numbers of the second kind given by
\begin{equation}
\frac{t}{\log_{\lambda} (1+t)}=\sum_{l=0}^{\infty}b_{l,\lambda} \frac{t^{l}}{l!}. \label{10e}
\end{equation}
From Theorem 3.3, \eqref{4e} and \eqref{9e}, we obtain the following expression
\begin{align*}
x^{n}&=\sum_{l=0}^{n}\frac{1}{p^{l-1}}S_{2}(n,l)b_{l,\lambda}\beta_{0,\lambda}^{Y}(x) \\
&\quad\quad +\sum_{k=1}^{n}\Big\{\frac{1}{k}\sum_{j=k-1}^{n-1}\frac{1}{p^{j}}
S_{1,\lambda}(j,k-1)\frac{1}{j!}\Delta^{j+1}0^{n} \Big\}\beta_{k,\lambda}^{Y}(x).
\end{align*}
(b) Let $Y$ be the binomial random variable with parameter $(m,p)$. Then the probability mass function of $Y$ is given by (see [29])
\begin{equation*}
p(i)=\binom{m}{i}p^{i}(1-p)^{m-i}, \quad i=0,1,\dots, m.
\end{equation*}
Here $m$ is a positive integer and $0< p \le 1$. Then, from [14], we have
\begin{align}
&f_{Y}(t)=\log\big(1+\frac{1}{p}(e^{\frac{t}{m}}-1)\big),\,\, f_{Y,\lambda}(t)=\log_{\lambda}\big(1+\frac{1}{p}(e^{\frac{t}{m}}-1) \big), \label{11e} \\
&S_{1}^{Y}(n,k)=\sum_{l=k}^{n}\sum_{i=l}^{n}\frac{1}{p^{l}}\frac{1}{m^{i}}S_{1}(l,k)S_{2}(i,l)S_1(n,i), \nonumber \\
&S_{1,\lambda}^{Y}(n,k)=\sum_{l=k}^{n}\sum_{i=l}^{n}\frac{1}{p^{l}}\frac{1}{m^{i}}S_{1,\lambda}(l,k)S_{2}(i,l)S_{1}(n,i). \nonumber
\end{align} \par
Let $x^{n}=\sum_{k=0}^{n}a_{k}B_{k}^{Y}(x)$.
Then, by proceeding as in (a), we can show that
\begin{equation}
\frac{t}{f_{Y}(t)}x^{n}=\frac{1}{m^{n-1}}\sum_{r=0}^{n}\sum_{l=0}^{r}\binom{n}{r}\frac{1}{p^{l-1}}S_{2}(r,l)b_{l}B_{n-r}(mx). \label{12e}
\end{equation}
From Theorem 3.1, \eqref{11e} and \eqref{12e}, we derive
\begin{align*}
x^{n}&=\frac{1}{m^{n}}\sum_{r=0}^{n}\sum_{l=0}^{r}\binom{n}{r}\frac{1}{n-r+1}\frac{1}{p^{l-1}}S_{2}(r,l)b_{l}\big(B_{n-r+1}(m)-B_{n-r+1}\big) B_{0}^{Y}(x) \\
&\quad +\sum_{k=1}^{n}\Big\{\frac{1}{k}\sum_{j=k-1}^{n-1}\sum_{l=k-1}^{j}\sum_{i=l}^{j}\frac{1}{p^{l}}\frac{1}{m^{i}}S_{1}(l,k-1)S_{2}(i,l)S_{1}(j,i)\frac{1}{j!}\Delta^{j+1}0^{n}\Big\}B_{k}^{Y}(x).
\end{align*} \par
Next, we let $x^{n}=\sum_{k=0}^{n}a_{k}\beta_{k,\lambda}^{Y}(x)$. Then we see that
\begin{equation}
\frac{t}{f_{Y,\lambda}(t)}x^{n}=\frac{1}{m^{n-1}}\sum_{r=0}^{n}\sum_{l=0}^{r}\binom{n}{r}\frac{1}{p^{l-1}}S_{2}(r,l)b_{l,\lambda}B_{n-r}(mx). \label{13e}
\end{equation}
From Theorem 3.3, \eqref{11e} and \eqref{13e}, we deduce
\begin{align*}
x^{n}&=\frac{1}{m^{n}}\sum_{r=0}^{n}\sum_{l=0}^{r}\binom{n}{r}\frac{1}{n-r+1}\frac{1}{p^{l-1}}S_{2}(r,l)b_{l,\lambda}\big(B_{n-r+1}(m)-B_{n-r+1}\big) \beta_{0,\lambda}^{Y}(x) \\
&\quad +\sum_{k=1}^{n}\Big\{\frac{1}{k}\sum_{j=k-1}^{n-1}\sum_{l=k-1}^{j}\sum_{i=l}^{j}\frac{1}{p^{l}}\frac{1}{m^{i}}S_{1,\lambda}(l,k-1)S_{2}(i,l)S_{1}(j,i)\frac{1}{j!}\Delta^{j+1}0^{n}\Big\}\beta_{k,\lambda}^{Y}(x).
\end{align*}
(c) Let $Y$ be the Poisson random variable with parameter $\alpha >0$. Then the probability mass function of $Y$ is given by (see [29])
\begin{equation*}
p(i)=e^{-\alpha}\frac{\alpha^{i}}{i!}, \quad i=0,1,2,\dots.
\end{equation*}
Then, from [14], we have
\begin{align}
&f_{Y}(t)=\log \big(1+\frac{t}{\alpha}\big), \,\, f_{Y,\lambda}(t)=\log_{\lambda}\big(1+\frac{t}{\alpha}\big), \label{14e} \\
&S_{1}^{Y}(n,k)=\sum_{l=k}^{n}\frac{1}{\alpha^{l}}S_{1}(l,k)S_{1}(n,l), \,\,
S_{1,\lambda}^{Y}(n,k)=\sum_{l=k}^{n}\frac{1}{\alpha^{l}}S_{1,\lambda}(l,k)S_{1}(n,l). \nonumber
\end{align} \par
Let $x^{n}=\sum_{k=0}^{n}a_{k}B_{k}^{Y}(x)$. Then we can show that
\begin{equation}
\frac{t}{f_{Y}(t)}x^{n}=\sum_{l=0}^{n}\binom{n}{l}\frac{1}{\alpha^{l-1}}b_{l}x^{n-l}.\label{15e}
\end{equation}
Then, from Theorem 3.1, \eqref{14e} and \eqref{15e}, we obtain
\begin{align*}
x^{n}&=\sum_{l=0}^{n}\binom{n}{l}\frac{1}{n-l+1}\frac{1}{\alpha^{l-1}}b_{l} B_{0}^{Y}(x) \\
&\quad +\sum_{k=1}^{n}\Big\{\frac{1}{k}\sum_{j=k-1}^{n-1}\sum_{l=k-1}^{j}\frac{1}{\alpha^{l}}S_{1}(l,k-1)S_{1}(j,l)\frac{1}{j!}\Delta^{j+1}0^{n}\Big\}B_{k}^{Y}(x).
\end{align*} \par
Next, we let $x^{n}=\sum_{k=0}^{n}a_{k}\beta_{k,\lambda}^{Y}(x)$. Then we see that
\begin{equation}
\frac{t}{f_{Y,\lambda}(t)}x^{n}=\sum_{l=0}^{n}\binom{n}{l}\frac{1}{\alpha^{l-1}}b_{l,\lambda}x^{n-l}.\label{16e}
\end{equation}
Then, from Theorem 3.3, \eqref{14e} and \eqref{16e}, we derive
\begin{align*}
x^{n}&=\sum_{l=0}^{n}\binom{n}{l}\frac{1}{n-l+1}\frac{1}{\alpha^{l-1}}b_{l,\lambda} \beta_{0,\lambda}^{Y}(x) \\
&\quad +\sum_{k=1}^{n}\Big\{\frac{1}{k}\sum_{j=k-1}^{n-1}\sum_{l=k-1}^{j}\frac{1}{\alpha^{l}}S_{1,\lambda}(l,k-1)S_{1}(j,l)\frac{1}{j!}\Delta^{j+1}0^{n}\Big\}\beta_{k,\lambda}^{Y}(x).
\end{align*}
(d) Let $Y$ be the geometric random variable with parameter $0 <p <1$. Then the probability mass function of $Y$ is given by (see [29])
\begin{equation*}
p(i)=(1-p)^{i-1}p, \quad i=1,2,\dots.
\end{equation*}
Then, from [14], we have
\begin{align}
&f_{Y}(t)=\log \Big(\frac{e^{t}}{p+(1-p)e^{t}}\Big),\,\, f_{Y,\lambda}(t)=\log_{\lambda} \Big(\frac{e^{t}}{p+(1-p)e^{t}}\Big), \label{17e} \\
&S_{1}^{Y}(n,k)=\sum_{j=k}^{n}\binom{n}{j}(n-1)_{n-j}p^{j}(p-1)^{n-j}S_{1}(j,k), \nonumber \\
&S_{1,\lambda}^{Y}(n,k)=\sum_{j=k}^{n}\binom{n}{j}(n-1)_{n-j}p^{j}(p-1)^{n-j}S_{1,\lambda}(j,k). \nonumber
\end{align} \par
Let $x^{n}=\sum_{k=0}^{n}a_{k}B_{k}^{Y}(x)$. Then we consider
\begin{equation}
\frac{t}{\log \Big(\frac{e^{t}}{p+(1-p)e^{t}}\Big)}e^{xt}=\frac{\frac{p(e^{t}-1)}{p+(1-p)e^{t}}}{\log\Big(1+\frac{p(e^{t}-1)}{p+(1-p)e^{t}}\Big)}\frac{t}{\frac{p(e^{t}-1)}{p+(1-p)e^{t}}}e^{xt}. \label{18e}
\end{equation}
Noting that $\frac{p(e^{t}-1)}{p+(1-p)e^{t}}=\frac{p}{1-p}\Big(1-\frac{1-\frac{p}{p-1}}{e^{t}-\frac{p}{p-1}}\Big)$ and from \eqref{35a}, the first factor on the right hand side of \eqref{18e} is
\begin{align}
\frac{\frac{p(e^{t}-1)}{p+(1-p)e^{t}}}{\log\Big(1+\frac{p(e^{t}-1)}{p+(1-p)e^{t}}\Big)}
&=\sum_{l=0}^{\infty}b_{l}\frac{1}{l!}\Big(\frac{p}{1-p}\Big)^{l}\Big(1-\frac{1-\frac{p}{p-1}}{e^{t}-\frac{p}{p-1}}\Big)^{l} \label{19e} \\
&=\sum_{l=0}^{\infty}b_{l}\frac{1}{l!}\Big(\frac{p}{1-p}\Big)^{l}\sum_{r=0}^{l}\binom{l}{r}(-1)^{r}\bigg(\frac{1-\frac{p}{p-1}}{e^{t}-\frac{p}{p-1}}\bigg)^{r} \nonumber \\
&=\sum_{l=0}^{\infty}b_{l}\frac{1}{l!}\Big(\frac{p}{1-p}\Big)^{l}\sum_{r=0}^{l}\binom{l}{r}(-1)^{r}\sum_{j=0}^{\infty}H_{j}^{(r)}\big(\frac{p}{p-1}\big)\frac{t^{j}}{j!} \nonumber \\
&=\sum_{j=0}^{\infty}\sum_{l=0}^{\infty}\sum_{r=0}^{l}(-1)^{r}\frac{1}{l!}\binom{l}{r}\Big(\frac{p}{1-p}\Big)^{l}b_{l}H_{j}^{(r)}\big(\frac{p}{p-1}\big)\frac{t^{j}}{j!}. \nonumber
\end{align}
Recalling that $B_{m}(x+1)-B_{m}(x)=mx^{m-1}$ (see \eqref{6c} with $Y=1$), the second factor of \eqref{18e} is
\begin{align}
\frac{t}{\frac{p(e^{t}-1)}{p+(1-p)e^{t}}}e^{xt}&=\frac{1}{p}\Big(p\frac{t}{e^{t}-1}e^{xt}+(1-p)\frac{t}{e^{t}-1}e^{(x+1)t}\Big) \label{20e} \\
&=\frac{1}{p}\sum_{m=0}^{\infty}\big(B_{m}(x+1)-pmx^{m-1} \big)\frac{t^{m}}{m!}.\nonumber
\end{align}
From \eqref{17e}-\eqref{20e}, we derive
\begin{align}
\frac{t}{f_{Y}(t)}x^{n}&=\frac{1}{p}\sum_{j=0}^{n}\sum_{l=0}^{\infty}\sum_{r=0}^{l}(-1)^{r}\frac{1}{l!}\binom{n}{j}\binom{l}{r}\Big(\frac{p}{1-p}\Big)^{l}b_{l} \label{21e} \\
&\quad\quad \times H_{j}^{(r)}\Big(\frac{p}{p-1}\Big)\big(B_{n-j}(x+1)-p(n-j)x^{n-j-1} \big). \nonumber
\end{align}
Before proceeding further, we note that
\begin{align}
\int_{0}^{1}&\big(B_{n-j}(x+1)-p(n-j)x^{n-j-1} \big) dx \label{22e} \\
&=\Big[\frac{1}{n-j+1}B_{n-j+1}(x+1)-px^{n-j} \Big]_{0}^{1} \nonumber \\
&=\frac{1}{n-j+1}\big(B_{n-j+1}(2)-B_{n-j+1}(1)\big)-p(1-\delta_{n,j}) \nonumber \\
&=1-p(1-\delta_{n,j}). \nonumber
\end{align}
So, from \eqref{21e} and \eqref{22e}, we get
\begin{align}
a_{0}&=\int_{0}^{1}\frac{t}{f_{Y}(t)}x^{n} dx=\frac{1}{p}\sum_{j=0}^{n}\sum_{l=0}^{\infty}\sum_{r=0}^{l}(-1)^{r}\frac{1}{l!}\binom{n}{j}\binom{l}{r}\Big(\frac{p}{1-p}\Big)^{l}b_{l} \label{23e} \\
&\quad\quad\quad\quad\quad\quad\quad\quad\quad \times H_{j}^{(r)}\Big(\frac{p}{p-1}\Big)\big(1-p(1-\delta_{n,j})\big). \nonumber
\end{align}
Now, from Theorem 3.1, \eqref{17e} and \eqref{23e}, we obtain
\begin{align*}
x^{n}&=\frac{1}{p}\sum_{j=0}^{n}\sum_{l=0}^{\infty}\sum_{r=0}^{l}(-1)^{r}\frac{1}{l!}\binom{n}{j}\binom{l}{r}\Big(\frac{p}{1-p}\Big)^{l}b_{l}\\
&\quad\quad \times H_{j}^{(r)}\Big(\frac{p}{p-1}\Big)\big(1-p(1-\delta_{n,j})\big)B_{0}^{Y}(x) \\
& \quad\quad +\sum_{k=1}^{n}\Big\{\frac{1}{k}\sum_{j=k-1}^{n-1}\sum_{l=k-1}^{j}\binom{j}{l}(j-1)_{j-l}p^{l}(p-1)^{j-l} \\
&\quad\quad \times S_{1}(l,k-1)\frac{1}{j!}\Delta^{j+1}0^{n}\Big\}B_{k}^{Y}(x).
\end{align*} \par
Next, we let $x^{n}=\sum_{k=0}^{n}a_{k}\beta_{k,\lambda}^{Y}(x)$. Then we see that
\begin{align}
a_{0}&=\int_{0}^{1}\frac{t}{f_{Y,\lambda}(t)}x^{n} dx=\frac{1}{p}\sum_{j=0}^{n}\sum_{l=0}^{\infty}\sum_{r=0}^{l}(-1)^{r}\frac{1}{l!}\binom{n}{j}\binom{l}{r}\Big(\frac{p}{1-p}\Big)^{l}b_{l,\lambda} \label{24e} \\
&\quad\quad\quad\quad\quad\quad\quad\quad\quad \times H_{j}^{(r)}\Big(\frac{p}{p-1}\Big)\big(1-p(1-\delta_{n,j})\big). \nonumber
\end{align}
Now, from Theorem 3.3, \eqref{17e} and \eqref{24e}, we get
\begin{align*}
x^{n}&=\frac{1}{p}\sum_{j=0}^{n}\sum_{l=0}^{\infty}\sum_{r=0}^{l}(-1)^{r}\frac{1}{l!}\binom{n}{j}\binom{l}{r}\Big(\frac{p}{1-p}\Big)^{l}b_{l,\lambda}\\
&\quad\quad \times H_{j}^{(r)}\Big(\frac{p}{p-1}\Big)\big(1-p(1-\delta_{n,j})\big)\beta_{0,\lambda}^{Y}(x) \\
& \quad\quad +\sum_{k=1}^{n}\Big\{\frac{1}{k}\sum_{j=k-1}^{n-1}\sum_{l=k-1}^{j}\binom{j}{l}(j-1)_{j-l}p^{l}(p-1)^{j-l} \\
&\quad\quad \times S_{1,\lambda}(l,k-1)\frac{1}{j!}\Delta^{j+1}0^{n}\Big\}\beta_{k,\lambda}^{Y}(x).
\end{align*}
(e) Let $Y$ be the exponential random variable with parameter $\alpha > 0$. Then the probability density function of $Y$ is given by (see [29])
\begin{equation*}
f(y)=\left\{\begin{array}{ccc}
\alpha e^{-\alpha y}, & \textrm{if \,\,$y \ge 0$,} \\
0, & \textrm{if\,\, $y<0$}.
\end{array}\right.
\end{equation*}
Then, from [14], we have
\begin{align}
&f_{Y}(t)=\alpha (1-e^{-t}),\,\,f_{Y,\lambda}(t)=\frac{1}{\lambda} \big( e^{\alpha \lambda (1-e^{-t})} -1 \big), \label{25e} \\
&S_{1}^{Y}(n,k)=(-1)^{n-k}\binom{n}{k}(n-1)_{n-k}\alpha^{k},\nonumber \\
&S_{1,\lambda}^{Y}(n,k)=\sum_{j=k}^{n}\binom{n}{j}(-1)^{n-j}(n-1)_{n-j}\alpha^{j} \lambda^{j-k}S_{2}(j,k). \nonumber
\end{align} \par
Let $x^{n}=\sum_{k=0}^{n}a_{k}B_{k}^{Y}(x)$. As $\frac{t}{\alpha (1-e^{-t})}e^{xt}=\frac{1}{\alpha}\sum_{n=0}^{\infty}(-1)^{n}B_{n}(-x)\frac{t^{n}}{n!}$, we have
\begin{equation}
\frac{t}{\alpha (1-e^{-t})}x^{n}=\frac{1}{\alpha}(-1)^{n}B_{n}(-x). \label{26e}
\end{equation}
Then, by using Lemma 5.1 and \eqref{26e}, we deduce
\begin{equation}
a_{0}=\int_{0}^{1}\frac{t}{f_{Y}(t)}x^{n}dx= \frac{1}{\alpha}(-1)^{n}\int_{0}^{1}B_{n}(-x) dx=\frac{1}{\alpha}. \label{27e}
\end{equation}
From Theorem 3.1, \eqref{25e} and \eqref{26e}, we derive
\begin{align*}
x^{n}&=\frac{1}{\alpha}B_{0}^{Y}(x)+\sum_{k=1}^{n}\Big\{\frac{1}{k}\sum_{j=k-1}^{n-1}(-1)^{j-k+1}\binom{j}{k-1} \\
&\quad\quad \times (j-1)_{j-k+1}\alpha^{k-1}\frac{1}{j!}\Delta^{j+1}0^{n}\Big\}B_{k}^{Y}(x).
\end{align*} \par
Next, we let $x^{n}=\sum_{k=0}^{n}a_{k}\beta_{k,\lambda}^{Y}(x)$. We note that
\begin{align}
\frac{t}{e^{\alpha \lambda (1-e^{-t})}-1}e^{xt}&=\frac{\alpha \lambda (1-e^{-t})}{e^{\alpha \lambda (1-e^{-t})}-1}\frac{t}{\alpha \lambda (1-e^{-t})}e^{xt} \label{28e} \\
&=\frac{1}{\alpha \lambda}\sum_{l=0}^{\infty}B_{l}(-\alpha \lambda)^{l}\frac{1}{l!}(e^{-t}-1)^{l}\sum_{m=0}^{\infty}B_{m}(-x)(-1)^{m}\frac{t^{m}}{m!} \nonumber \\
&=\frac{1}{\alpha \lambda}\sum_{l=0}^{\infty}B_{l}(-\alpha \lambda)^{l}\sum_{r=l}^{\infty}S_{2}(r,l)(-1)^{r}\frac{t^{r}}{r!}\sum_{m=0}^{\infty}B_{m}(-x)(-1)^{m}\frac{t^{m}}{m!} \nonumber \\
&=\frac{1}{\alpha \lambda}\sum_{n=0}^{\infty}\Big\{\sum_{r=0}^{n}\sum_{l=0}^{r}\binom{n}{r}S_{2}(r,l)(-1)^{l+n}(\alpha \lambda)^{l} B_{l}B_{n-r}(-x)\Big\}\frac{t^{n}}{n!}. \nonumber
\end{align}
Thus we have shown that
\begin{equation}
\frac{t}{e^{\alpha \lambda (1-e^{-t})}-1}x^{n}=\frac{1}{\alpha \lambda}\sum_{r=0}^{n}\sum_{l=0}^{r}\binom{n}{r}S_{2}(r,l)(-1)^{l+n}(\alpha \lambda)^{l} B_{l}B_{n-r}(-x). \label{29e}
\end{equation}
By using Lemma 5.1, \eqref{25e} and \eqref{29e}, we get
\begin{align}
a_{0}&=\int_{0}^{1}\frac{t}{f_{Y,\lambda}(t)}x^{n}dx=\lambda\int_{0}^{1}\frac{t}{e^{\alpha \lambda (1-e^{-t})}-1}x^{n} dx \label{30e} \\
&=\frac{1}{\alpha}\sum_{r=0}^{n}\sum_{l=0}^{r}\binom{n}{r}S_{2}(r,l)(-1)^{l-r}(\alpha \lambda)^{l} B_{l}. \nonumber
\end{align}
Now, from Theorem 3.3, \eqref{25e} and \eqref{30e}, we obtain
\begin{align*}
x^{n}&=\frac{1}{\alpha}\sum_{r=0}^{n}\sum_{l=0}^{r}\binom{n}{r}S_{2}(r,l)(-1)^{l-r}(\alpha \lambda)^{l} B_{l}\beta_{0,\lambda}^{Y}(x) \\
&\quad +\sum_{k=1}^{n}\Big\{\frac{1}{k}\sum_{j=k-1}^{n-1}\sum_{m=k-1}^{j}\binom{j}{m}(-1)^{j-m}(j-1)_{j-m}  \\
& \quad \quad \quad \times \alpha^{m}\lambda^{m-k+1}S_{2}(m,k-1)\frac{1}{j!}\Delta^{j+1}0^{n} \Big\}\beta_{k,\lambda}^{Y}(x).
\end{align*}
(f) Let $Y$ be the gamma random variable with parameters $\alpha, \beta >0$. Then the probability density function of $Y$ is given by (see [29])
\begin{equation*}
f(y)=\left\{\begin{array}{ccc}
\beta e^{-\beta y}\frac{(\beta y)^{\alpha-1}}{\Gamma(\alpha)}, & \textrm{if\,\, $y\ge 0$,}\\
0, & \textrm{if \,\,$y<0$},
\end{array}\right.
\end{equation*}
Then, from [14], we have
\begin{align}
&f_{Y}(t)=\beta\big(1-e^{-\frac{t}{\alpha}} \big),\,\,f_{Y, \lambda}(t)=\frac{1}{\lambda}\Big(e^{\lambda \beta(1-e^{-\frac{t}{\alpha}})}-1 \Big), \label{31e} \\
&S_{1}^{Y}(n,k)=\frac{1}{k!}\beta^{k}\big(-\frac{1}{\alpha} \big)^{n}\sum_{l=0}^{k}\binom{k}{l}(-1)^{l}\big(l+(n-1)\alpha \big)_{n,\alpha}, \nonumber \\
&S_{1,\lambda}^{Y}(n,k)=\sum_{l=k}^{n}\sum_{r=0}^{l}S_{2}(l,k)\binom{l}{r}\frac{1}{l!}(-1)^{n-r}\lambda^{l-k}\beta^{l}\frac{1}{\alpha^{n}}\big(r+(n-1)\alpha \big)_{n,\alpha}. \nonumber
\end{align} \par
Let $x^{n}=\sum_{k=0}^{n}a_{k}B_{k}^{Y}(x)$. Observe first that
\begin{equation}
\frac{t}{\beta (1-e^{-\frac{t}{\alpha}})}e^{xt}=\frac{\alpha}{\beta}\sum_{n=0}^{\infty}B_{n}(-\alpha x)\Big(-\frac{1}{\alpha}\Big)^{n}\frac{t^{n}}{n!}. \label{32e}
\end{equation}
Thus, from \eqref{31e} and \eqref{32e}, we obtain
\begin{equation}
\frac{t}{f_{Y}(t)}x^{n}=\frac{\alpha}{\beta}\Big(-\frac{1}{\alpha}\Big)^{n}B_{n}(-\alpha x), \label{33e}
\end{equation}
which yields that
\begin{equation}
a_{0}=\int_{0}^{1}\frac{t}{f_{Y}(t)}x^{n} dx=\frac{\alpha}{\beta}\Big(-\frac{1}{\alpha} \Big)^{n+1}\frac{1}{n+1}\big(B_{n+1}(-\alpha)-B_{n+1} \big). \label{34e}
\end{equation}
From Theorem 3.1, \eqref{31e} and \eqref{34e}, we get
\begin{align*}
x^{n}&=\frac{\alpha}{\beta}\Big(-\frac{1}{\alpha} \Big)^{n+1}\frac{1}{n+1}\big(B_{n+1}(-\alpha)-B_{n+1} \big)B_{0}^{Y}(x) \\
&\quad +\sum_{k=1}^{n}\Big(\frac{1}{k!}\beta^{k-1}\sum_{j=k-1}^{n-1}\sum_{l=0}^{k-1}\Big(-\frac{1}{\alpha}\Big)^{j}\binom{k-1}{l}(-1)^{l} \\
&\quad \times \big(l+(j-1)\alpha \big)_{j,\alpha}\frac{1}{j!}\Delta^{j+1}0^{n} \Big) B_{k}^{Y}(x).
\end{align*} \par
Next, we let $x^{n}=\sum_{k=0}^{n}a_{k}\beta_{k,\lambda}^{Y}(x)$. We note that (see \eqref{32e})
\begin{align}
&\frac{t}{\frac{1}{\lambda}(e^{\lambda \beta(1-e^{-\frac{t}{\alpha}})}-1)}e^{xt}=
\frac{\lambda \beta(1-e^{-\frac{t}{\alpha}})}{e^{\lambda \beta(1-e^{-\frac{t}{\alpha}})}-1}\frac{t}{\beta (1-e^{-\frac{t}{\alpha}})}e^{xt} \label{35} \\
&=\frac{\alpha}{\beta}\sum_{l=0}^{\infty}B_{l}(-\lambda \beta)^{l}\frac{1}{l!}\big(e^{-\frac{t}{\alpha}}-1 \big)^{l}\sum_{m=0}^{\infty}B_{m}(-\alpha x)\Big(-\frac{1}{\alpha}\Big)^{m}\frac{t^{m}}{m!} \nonumber \\
&=\frac{\alpha}{\beta}\sum_{l=0}^{\infty}B_{l}(-\lambda \beta)^{l}\sum_{r=l}^{\infty}S_{2}(r,l)\Big(-\frac{1}{\alpha}\Big)^{r}\frac{t^{r}}{r!}\sum_{m=0}^{\infty}B_{m}(-\alpha x)\Big(-\frac{1}{\alpha}\Big)^{m}\frac{t^{m}}{m!} \nonumber \\
&=\frac{\alpha}{\beta}\sum_{n=0}^{\infty}\Big\{\sum_{r=0}^{n}\sum_{l=0}^{r}\binom{n}{r}(-\lambda \beta)^{l}\Big(-\frac{1}{\alpha}\Big)^{n}S_{2}(r,l)B_{l}B_{n-r}(-\alpha x) \Big\}\frac{t^{n}}{n!}, \nonumber
\end{align}
from which we deduce
\begin{equation}
\frac{t}{f_{Y,\lambda}(t)}x^{n}=\frac{\alpha}{\beta}\sum_{r=0}^{n}\sum_{l=0}^{r}\binom{n}{r}(-\lambda \beta)^{l}\Big(-\frac{1}{\alpha}\Big)^{n}S_{2}(r,l)B_{l}B_{n-r}(-\alpha x). \label{36e}
\end{equation}
Now, from Theorem 3.3, \eqref{31e} and \eqref{36e}, we derive
\begin{align*}
x^{n}&=\Big\{\frac{\alpha}{\beta}\sum_{r=0}^{n}\sum_{l=0}^{r}\frac{1}{n-r+1}\binom{n}{r}(-\lambda \beta)^{l}\Big(-\frac{1}{\alpha}\Big)^{n+1} \\
&\quad\quad \times S_{2}(r,l)B_{l}\big(B_{n_r+1}(-\alpha)-B_{n-r+1} \Big\}\beta_{0,\lambda}^{Y}(x) \\
&\quad +\sum_{k=1}^{n}\Big\{\frac{1}{k}\sum_{j=k-1}^{n-1}\sum_{l=k-1}^{j}\sum_{r=0}^{l}S_{2}(l,k-1)\binom{l}{r}\frac{1}{l!}(-1)^{j-r}\lambda^{l-k+1} \\
&\quad\quad \times \beta^{l}\frac{1}{\alpha^{j}}\big(r+(j-1)\alpha\big)_{j,\alpha}\frac{1}{j!}\Delta^{j+1}0^{n} \Big\}\beta_{k,\lambda}^{Y}(x).
\end{align*}

\section{Conclusion}
This paper explored the representation of arbitrary polynomials as linear combinations of probabilisitc Bernoulli polynomials associated with $Y$, $B_{k}^{Y}(x)$, and probabilistic degenerate Bernoulli polynomials associated with $Y$, $\beta_{k,\lambda}^{Y,(r)}(x)$, and their higher-order counterparts, $B_{k}^{Y,(r)}(x)$ and $\beta_{k,\lambda}^{Y,(r)}(x)$. We derived explicit coefficients for these linear combinations, expressed in terms of probabilistic Stirling numbers of the first kind associated with $Y$, $S_{1}^{Y}(n,k)$, and probabilistic degenerate Stirling numbers of the first kind associated with $Y$ $S_{1,\lambda}^{Y}(n,k)$. To demonstrate our findings, we provided concrete examples by expressing $x^{n}$ as linear combinations of these polynomials for various discrete and continuous random variables $Y$, utilizing established results for $f_{Y}(t),\,\,f_{Y,\lambda}(t),\,\,S_{1}^{Y}(n,k),\,\,and\,\, S_{1,\lambda}^{Y}(n,k)$, where the compositional inverses of $f_{Y}(t)$ and $f_{Y,\lambda}(t)$ are respectively given by
\begin{equation*}
\bar{f}_{Y}(t)=\log E[e^{Y t}], \quad \bar{f}_{Y,\lambda}(t)=\log E[e_{\lambda}^{Y}(t)].
\end{equation*}

\noindent{\bf Conflict of interest}\\
 The authors declare that they have no Conflict of interest.

 %\noindent{\bf Acknowledgement}\\

%----------------------------------

\end{document}